\documentclass[12pt,reqno]{amsart}

\usepackage{amssymb}
\usepackage{bbm}
\usepackage{mathrsfs}
\usepackage{xfrac}
\usepackage{anyfontsize}
\usepackage{xypic}
\usepackage{float}

\newtheorem{facs}{Facts}[section]
\newtheorem{theo}[facs]{Theorem}
\newtheorem{prop}[facs]{Proposition}
\newtheorem{lem}[facs]{Lemma}
\newtheorem{coro}[facs]{Corollary}

\newtheorem*{theoo}{Theorem}

\theoremstyle{definition}
\newtheorem{algo}[facs]{Algorithm}
\newtheorem{defi}[facs]{Definition}
\newtheorem{ttt}[facs]{}

\newtheorem*{algoo}{Algorithm}
\newtheorem*{nota}{Notation}
\newtheorem*{notao}{Terminology and Notation}
\newtheorem*{inp}{Input}

\theoremstyle{remark}
\newtheorem{rem}[facs]{Remark}
\newtheorem{rems}[facs]{Remarks}
\newtheorem{ex}[facs]{Example}
\newtheorem{exs}[facs]{Examples}

\makeatletter
\def\hsmash{\relax 
  \ifmmode\def\next{\mathpalette\mathhsm@sh}\else\let\next\makehsm@sh
  \fi\next}
\def\makehsm@sh#1{\setbox\z@\hbox{#1}\finhsm@sh}
\def\mathhsm@sh#1#2{\setbox\z@\hbox{$\m@th#1{#2}$}\finhsm@sh}
\def\finhsm@sh{\wd\z@\z@ \box\z@}
\makeatother

\newcommand{\brr}{, }
\newcommand{\br}{ }

\newcommand{\Gal}{\mathop{\text{\rm Gal}}\nolimits}
\newcommand{\Div}{\mathop{\text{\rm Div}}\nolimits}
\newcommand{\Pic}{\mathop{\text{\rm Pic}}\nolimits}
\newcommand{\Frob}{\mathop{\text{\rm Frob}}\nolimits}
\newcommand{\Spec}{\mathop{\text{\rm Spec}}\nolimits}
\newcommand{\Br}{\mathop{\text{\rm Br}}\nolimits}
\newcommand{\Hom}{\mathop{\text{\rm Hom}}\nolimits}
\newcommand{\GL}{\mathop{\text{\rm GL}}\nolimits}
\newcommand{\End}{\mathop{\text{\rm End}}\nolimits}
\newcommand{\Aut}{\mathop{\text{\rm Aut}}\nolimits}
\newcommand{\rk}{\mathop{\text{\rm rk}}\nolimits}
\newcommand{\ord}{\mathop{\text{\rm ord}}\nolimits}
\newcommand{\Alt}{\mathop{\text{\rm Alt}}\nolimits}
\newcommand{\Sym}{\mathop{\text{\rm Sym}}\nolimits}

\newcommand{\Mat}{\mathop{\text{\rm M}}\nolimits}

\newcommand{\T}{\mathop{\text{\rm T}}\nolimits}
\renewcommand{\O}{\mathop{\text{\rm O}}\nolimits}
\renewcommand{\P}{\mathop{\text{\rm P}}\nolimits}
\newcommand{\V}{\mathop{\text{\rm V}}\nolimits}
\newcommand{\Tr}{\mathop{\text{\rm Tr}}\nolimits}
\newcommand{\im}{\mathop{\text{\rm im}}\nolimits}

\newcommand{\et}{\text{\rm \'et}}

\newcommand{\id}{\text{\rm id}}
\newcommand{\even}{\text{\rm even}}
\newcommand{\alg}{\text{\rm alg}}

\newcommand{\Pb}{\mathop{\text{\bf P}}\nolimits}

\newcommand{\bbC}{{\mathbbm C}}
\newcommand{\bbF}{{\mathbbm F}}
\newcommand{\bbG}{{\mathbbm G}}
\newcommand{\bbN}{{\mathbbm N}}
\newcommand{\bbP}{{\mathbbm P}}
\newcommand{\bbQ}{{\mathbbm Q}}
\newcommand{\bbZ}{{\mathbbm Z}}

\newcommand{\calB}{{\mathscr{B}}}
\newcommand{\calL}{{\mathscr{L}}}
\newcommand{\calO}{{\mathscr{O}}}

\renewcommand{\textonehalf}{{\sfrac{1}{2}}}

\newcommand{\wpsi}{\rlap{\raisebox{-0.5ex}{$\:\!\widetilde{\phantom\psi}$}}\psi}

\renewcommand{\atop}[2]{\genfrac{}{}{0pt}{}{#1}{#2}}

\newcommand{\pmodulo}[1]{\nobreak\mkern8mu
 (\textup{mod}\,\,#1)}

\newcounter{ABC}
\newenvironment{ABC}{\begin{list}{\rm \Alph{ABC}) }%
{\usecounter{ABC} \leftmargin=0.0pt \labelsep=0.0pt %
\listparindent=0.0pt \labelwidth=0.0pt \parsep=\smallskipamount%
 \itemsep=0.0pt \topsep=0.0pt \partopsep=\smallskipamount}}{\end{list}}

\newcounter{abc}
\newenvironment{abc}{\begin{list}{\rm \alph{abc}) }%
{\usecounter{abc} \leftmargin=0.0pt \labelsep=0.0pt %
\listparindent=0.0pt \labelwidth=0.0pt \parsep=\smallskipamount%
 \itemsep=0.0pt \topsep=0.0pt \partopsep=\smallskipamount}}{\end{list}}

\newcounter{iii}
\newenvironment{iii}{\begin{list}{\rm \roman{iii}) }%
{\usecounter{iii} \leftmargin=0.0pt \labelsep=0.0pt %
\listparindent=0.0pt \labelwidth=0.0pt \parsep=\smallskipamount%
 \itemsep=0.0pt \topsep=0.0pt \partopsep=\smallskipamount}}{\end{list}}

\def\rightend#1#2{{%
 \leavevmode\nobreak\hskip .5em plus 1fil
 \penalty600 \hskip 0pt plus -1filll
 \vadjust{}\nobreak\hskip 0pt plus 1filll%
 #1\parfillskip=#2\relax \par}}

\def\eop{\ifmmode\rule[-22pt]{0pt}{1pt}\ifinner\tag*{$\square$}\else\eqno{\square}\fi\else\rightend{$\square$}{0pt}\fi}

\makeatletter
\newcommand{\leqnomode}{\tagsleft@true}
\newcommand{\reqnomode}{\tagsleft@false}
\makeatother

\title[2-adic point counting on
$K3$
surfaces]{2-adic point counting on
\boldmath{$K3$}
surfaces}

\begin{document}

\author{Andreas-Stephan Elsenhans}

\address{Institut f\"ur Mathematik\\ Universit\"at W\"urzburg\\ Emil-Fischer-Stra\ss e 30\\ D-97074 W\"urzburg\\ Germany}
\email{stephan.elsenhans@mathematik.uni-wuerzburg.de}
\urladdr{\;https://www.mathematik.uni-wuerzburg.de/computeralgebra/team/elsenhans-step\discretionary{}{}{}han-\discretionary{}{}{}prof-dr/}

\author{J\"org Jahnel}

\address{\mbox{Department \!Mathematik\\ \!Univ.\ \!Siegen\\ \!Walter-Flex-Str.\ \!3\\ D-57068 \!Siegen\\ \!Germany}}
\email{jahnel@mathematik.uni-siegen.de}
\urladdr{http://www.uni-math.gwdg.de/jahnel}


\date{June~28,~2022}



\begin{abstract}
This article reports on an approach to point counting on algebraic varieties over finite fields that is based on a detailed investigation of the
$2$-adic
orthogonal group. Combining the new approach with a
$p$-adic
method, we count the number of points on some
$K3$~surfaces
over the field
$\bbF_{\!p}$,
for all primes
$p < 10^8$.
\end{abstract}

\maketitle

\section{Introduction}

Counting points on algebraic varieties over finite fields is an important problem in algorithmic arithmetic geometry. When the Betti numbers of a variety are known, one has strong estimates on the number of points using \'etale cohomology.
For~example, in the case of a
$K3$~surface
$S$
that is projective over a finite field
$\bbF_{\!q}$,
the Lefschetz trace formula \cite[Rapport, Th\'eor\`eme~3.2]{SGA412} reads
\begin{equation}
\label{Lefschetz}
\#S(\bbF_{\!q}) = q^2 + \Tr(\Frob\colon H^2_\et(S_{\overline\bbF_{\!q}}, \bbZ_2(1)) \to H^2_\et(S_{\overline\bbF_{\!q}}, \bbZ_2(1))) q + 1 \,.
\end{equation}
Moreover,~according to the Weil conjectures proven by P.\ Deligne~\cite[Th\'e\-o\-r\`eme~1.6]{De74}, all eigenvalues of Frobenius are algebraic numbers of absolute
value~$1$.
As~$K3$
surfaces have
$\smash{\rk H^2_\et(S_{\overline\bbF_{\!q}}, \bbZ_2(1)) = 22}$
and the hyperplane section causes one eigenvalue to
be~$1$,
the inequality
\begin{equation}
\label{Weil_bound}
|\#S(\bbF_{\!q}) - (q^2 + q + 1)| \leq 21q
\end{equation}
results, which is of exactly the same form as the one formulated by P.\ Deligne for hypersurfaces \cite[Th\'e\-o\-r\`eme~8.1]{De74}. In~particular, one sees that it is sufficient to determine
$\#S(\bbF_{\!q})$
modulo some auxiliary integer that is larger
than~$42q$.

Nowadays,~the
$p$-adic
methods, as developed by K.\ Kedlaya, D.\ Harvey, and others, are frequently used for point counting, cf.\ \cite{Ke,Ha,HS,EJ16}. They~determine the number
$\#S(\bbF_{\!q})$,
for~$q$
a power of the prime
number~$p$,
by actually computing
\mbox{$(\#S(\bbF_{\!q}) \bmod p^j)$},
for a suitable value of the
exponent~$j$.

In this note, we are interested only in the number
$\#S(\bbF_{\!p})$
of points over the prime field
$\bbF_{\!p}$.
Then~the estimate~(\ref{Weil_bound}) shows that, in most cases,
$(\#S(\bbF_{\!p}) \bmod p^2)$
carries enough information.
However, for most of the primes, even the modulus
$p^2$
is by far larger than necessary. Moreover, the
$p$-adic
methods are faster, at least by a factor
of~$10$,
when working only
modulo~$p$
and not
modulo~$p^2$,
cf.\ Remark~\ref{factor10}. It~is thus worth trying to use a
$p$-adic
method just for counting
modulo~$p$
and to combine the result with the point count modulo some other small integer, which has to be obtained in a different~way.

\subsubsection*{
$l$-adic
point counting in general -- Explicitisation of \'etale cohomology}
Let~$l$
be a prime that is not the characteristic of the base~field. Then, the essence of an
\mbox{$l$-adic}
point counting algorithm is to make the
\mbox{$l$-adic}
cohomology
\mbox{$\bbZ_l$-module}
$H^i_\et(S_{\overline\bbF_{\!p}}, \bbZ_l)$
explicit, including the action
of~$\Frob$,
for the varieties of type~considered.

A~famous example is R.~Schoof's algorithm~\cite{Sch} for elliptic curves and its generalisation to abelian varieties~\cite{Pi}. Here,~the
\mbox{$l$-adic}
cohomology may be explicitly described using torsion~points.
Another well-known example works for del Pezzo surfaces. Here, formula~(\ref{Lefschetz}) holds, as~well. Moreover,~one has
$\smash{\T(S_{\overline\bbF_{\!p}}, \bbZ_l) = 0}$.
In~other words, 
$\smash{H^2_\et(S_{\overline\bbF_{\!p}}, \bbZ_l(1))}$
and the action
of~$\Frob$
upon it can be made explicit in terms of the geometric Picard group, cf.\ formula~(\ref{Lefschetz2}). This~essentially breaks down to the computation of the exceptional curves, cf.~\cite[\S2.5]{EJ20b}.

\subsubsection*{A\/
$2$-adic
point counting method for
$K3$~surfaces}
Following~these ideas, the algorithm we describe below relies on an explicitisation of
$\smash{\T(S_{\overline\bbF_{\!p}}, \bbZ_2) / 4 \T(S_{\overline\bbF_{\!p}}, \bbZ_2)}$
for a particular family~of
$K3$~surfaces.
This~turns out to suffice for point counting
modulo~$16$.

We~assume that we are given a
$K3$~surface
that is presented as a double cover
of~$\smash{\Pb^2_{\bbF_{\!p}}}$,
branched over six
$\bbF_{\!p}$-rational
lines. This~assumption is certainly more restrictive than necessary, but coincides with the generality, for which the algorithm is currently~implemented. It~coincides, too, with the generality, in which we describe the algorithm in Section~\ref{Algorithm}. We~indicate in Section~\ref{general} how to treat a slightly more general case and discuss the possibility of further generalisations in Section~\ref{very_general}.

The~very first step is to choose a lift to a flat
\mbox{$\bbZ$-scheme}~$S$
containing the given surface as the special
fibre~$\smash{S_{\bbF_{\!p}}}$.
Such~a lift exists for any
$K3$~surface,
at least as long as
$p \geq 5$,
as follows from \cite[Corol\-lary~2.3]{Og}, together with \cite[Theorem~1]{Ch}, cf.~\cite[Remarque~1.9]{De81}. In~our situation, one simply needs to lift the coefficients of the linear forms defining the lines from
$\bbF_{\!p}$
to~$\bbZ$.
One~then finds a double cover
$S$
of~$\smash{\Pb^2_\bbZ}$,
branched over the union of six lines, each of which is defined
over~$\Spec\bbZ$.
Then,~for every good prime
$l \neq 2$
of~$S$,
in particular for
$l = p$,
one has an isomorphism of
\mbox{$\bbZ_2$-modules}
$$H^2_\et(S_{\overline\bbF_{\!l}}, \bbZ_2(1)) \cong H^2_\et(S_{\overline\bbQ}, \bbZ_2(1)) \,,$$
the action of
$\Frob$
on the left agreeing with that
of~$\Frob_l$
on the right \cite[Expos\'e XVI, Corollaire~2.3]{SGA4}.

The~assumptions and preparations made up to here have two consequences, as explained in i) and~ii), below.

\begin{iii}
\item
One has that
$\smash{\rk\Pic(S_{\overline\bbQ}) \geq 16}$.
In~fact, a sublattice of
$\smash{\Pic(S_{\overline\bbQ})}$
of
rank~$16$
is explicitly known, which is a trivial
$\smash{\Gal(\overline\bbQ/\bbQ)}$-module.
Thus,~there is an improvement over~(\ref{Weil_bound}) implying that it is sufficient to count
$\smash{(\#S(\bbF_{\!p}) \bmod 16p)}$,
cf.\ Paragraph~\ref{Weil_bound_impr}. Assuming~that
$\smash{(\#S(\bbF_{\!p}) \bmod p)}$
is known, only the information about
$\smash{(\#S(\bbF_{\!p}) \bmod 16)}$
is~missing.

We~assume that
$\smash{\Pic(S_{\overline\bbQ})}$
is a trivial
$\smash{\Gal(\overline\bbQ/\bbQ)}$-module
in the case
$\smash{\rk\Pic(S_{\overline\bbQ}) > 16}$,
as~well. Then~in order to count
$\smash{(\#S(\bbF_{\!p}) \bmod 16)}$,
it suffices to~determine
$$\smash{(\Tr(\Frob_p\colon \T(S_{\overline\bbQ}, \bbZ_2) \to \T(S_{\overline\bbQ}, \bbZ_2)) \bmod 16) \,,}$$
for
$\smash{\T(S_{\overline\bbQ}, \bbZ_2) \subset H^2_\et(S_{\overline\bbQ}, \bbZ_2(1))}$
the transcendental lattice. Cf.\ Definition~\ref{transc_lattice} and formula~(\ref{Lefschetz3}).
\item
The
$\smash{\Gal(\overline\bbQ/\bbQ)}$-module
$\smash{\T(S_{\overline\bbQ}, \bbZ_2)/2\T(S_{\overline\bbQ}, \bbZ_2)}$
is~trivial. Indeed,~one has a canonical
$\smash{\Gal(\overline\bbQ/\bbQ})$-equivariant
isomorphism
$\smash{\Br(S_{\overline\bbQ})_2 \cong \Hom(\T(S_{\overline\bbQ}, \bbZ_2), \bbZ/2\bbZ)}$,
cf.\ Theorem~\ref{prop_Br2}. Moreover,~the
\mbox{$2$-torsion}
of the geometric Brauer group is well understood for double covers, thanks to the work of A.\,N.\ Skorobogatov~\cite[Theorem~1.1]{Sk}.
\end{iii}

Here,~an important observation comes into~play. The~action
of~$\smash{\Gal(\overline\bbQ/\bbQ)}$ on
$\smash{\T(S_{\overline\bbQ}, \bbZ_2)}$
takes place via maps being orthogonal with respect to the cup product pairing (Cf.\ Section~\ref{EtCohBr}). And~for the orthogonal group, a remarkable phenomenon of
\mbox{$2$-adic}
overdetermination occurs. We~discuss this in detail in Section~\ref{overdet}, which is, from the technical point of~view, the main part of this article. In~fact, the following is~true.

\begin{theoo}
Let\/~$n\in\bbN$.
With~respect to a non-degenerate, symmetric bilinear form
on\/~$\bbQ_2^n$,
let\/
$U_1, U_2 \in \Mat_{n \times n}(\bbZ_2)$
be orthogonal matrices such~that
$U_1 \equiv U_2 \pmodulo 4$.
If\/
$U_1 \equiv E_n \pmodulo 2$
then\/~$\Tr(U_1) \equiv \Tr(U_2) \pmodulo {16}$.
\end{theoo}

In~particular, for an orthogonal endomorphism\/
$a$
of\/~$\smash{\T(S_{\overline\bbQ}, \bbZ_2)}$
such that\/
\mbox{$a \equiv \id \pmod 2$},
the reduction\/
$\smash{(a \bmod 4) \in \End(\T(S_{\overline\bbQ}, \bbZ_2)/4\T(S_{\overline\bbQ}, \bbZ_2))}$
completely determines\/
$\smash{(\Tr(a\colon \T(S_{\overline\bbQ}, \bbZ_2) \to \T(S_{\overline\bbQ}, \bbZ_2)) \bmod 16)}$.

Moreover,
$\{A \in \GL_n(\bbZ / 4 \bbZ) \mid\! A \equiv E_n\!\! \pmodulo 2\}$
is an elementary abelian
\mbox{$2$-group}.
Therefore,~the splitting field
$K$
of the
$\smash{\Gal(\overline\bbQ/\bbQ)}$-action
on
$\smash{\T(S_{\overline\bbQ}, \bbZ_2)/4\T(S_{\overline\bbQ}, \bbZ_2)}$
is an abelian extension
of~$\bbQ$
of
exponent~$2$.
Furthermore,~$K$
is unramified at every good prime
$l \neq 2$
of~$S$,
so that one has
$\smash{K \subseteq \bbQ(\sqrt{-1}, \sqrt{2}, \sqrt{\mathstrut d} \mid d {\rm \;bad\;prime\;of\; }S)}$.
Thus,~in order to determine
$\smash{(\Tr(\Frob_p\colon \T(S_{\overline\bbQ}, \bbZ_2) \to \T(S_{\overline\bbQ}, \bbZ_2)) \bmod 16)}$,
it~suffices

\begin{iii}
\item
to look for a small good prime
$l$
such that the action of
$\Frob_l$
agrees with that
of~$\Frob_p$
on
$\smash{K = \bbQ(\sqrt{-1}, \sqrt{2}, \sqrt{d} \mid d {\rm \;bad\;prime\;of\; }S)}$
and
\item
to count
$\#S(\bbF_{\!l})$
by another method, either naive or
\mbox{$p$-adic},
and to deduce
$\smash{\Tr(\Frob\colon \T(S_{\overline\bbF_{\!l}}, \bbZ_2) \to \T(S_{\overline\bbF_{\!l}}, \bbZ_2))}$
from that value.
\end{iii}
A~situation, where this approach is particularly efficient, is when a
\mbox{$\bbZ$-scheme}
$S$
is given, the points on many special fibres
$\smash{S_{\bbF_{\!p}}}$
of which are to be~counted.

\subsubsection*{Practical experiments}
We applied the
\mbox{$2$-adic}
method, as described in combination with a
\mbox{$p$-adic}
method, to count the
$\bbF_{\!p}$-rational
points on some
$K3$~surfaces,
for all primes
$p$
up
to~$10^8$.
All~computations were done using {\tt magma}~\cite{BCP}.

\subsubsection*{Checks for correctness of the implementation}
For~each of the surfaces, we compared the output with the result of a completely naive algorithm, for all primes
$p<1000$,
and with that of a
\mbox{$p$-adic}
method, counting
modulo~$p^2$,
for all primes
$p<100\,000$.

There is a further check as follows. The~estimate~(\ref{DelWeil}) allows
for~$\#S(\bbF_{\!p})$
an interval of
length~$12p$.
Thus, knowing
$(\#S(\bbF_{\!p}) \bmod p)$,
there are only twelve (or
perhaps~$13$)
options
for~$\#S(\bbF_{\!p})$, which means that not all residues
modulo~$16$
are permissible. At~least a ``random'' bug would certainly produce such residue classes from time to~time.

\subsubsection*{Code}
On~the web pages of either author, related to the project described in this article, the following code is publicly available.

\begin{iii}
\item
Naive point counting, as used for initialisation and checks.
\item
A Harvey style
\mbox{$p$-adic}
point counting method in
\mbox{$p$-adic}
precision~$1$
with remainder tree, for surfaces of the shape
$w^2 = xyz f_3$.
It runs through all prime numbers
$p \leq 10^6$
in a few~minutes.
\item
The~initialisation following the more efficient approach, as described in Algorithm~\ref{2adic_count}.A') and~A'').
\item
Point counting
modulo~$16p$,
as described in Algorithm~\ref{2adic_count}.B).
\end{iii}

\begin{notao}
\begin{iii}
\item
For~$\calO$
a commutative ring
with~$1$
and
$n\in\bbN$,
we write
$\GL_n(\calO) := \{A\in\Mat_{n\times n}(\calO) \mid \det A {\rm ~is~a~unit~in~} \calO\}$.

We let
$E_n$
be the
$n\times n$
identity~matrix. For~any matrix
$A$,
we denote by
$A^\top$
the transpose in the usual~sense. Note~that
$A^\top$
is usually {\em not\/} the adjoint matrix in situations when the latter is~defined.

When
$\varphi\colon \Gamma \to \Gamma'$
is an
\mbox{$\calO$-linear}
map between free
$\calO$-modules
of finite rank with bases
$\calB$
and~$\calB'$,
respectively, then we denote by
$\Mat^{\calB'}_\calB(\varphi)$
the matrix
of~$\varphi$
with respect to
$\calB$
and~$\calB'$.
\item
By an
{\em $\calO$-lattice,}
we mean a free
\mbox{$\calO$-module}
$\Gamma$
of finite rank equipped with a non-degenerate, symmetric
\mbox{$\calO$-bilinear}
form
$b\colon \Gamma\times\Gamma \to \calO$.
When~there seems to be no danger of confusion, we simply write
$\Gamma$
instead
of~$(\Gamma, b)$.

An
$\calO$-lattice
$(\Gamma, b)$
is called {\em regular\/} if
$b$
provides an isomorphism
$\Gamma \to \Hom_\calO(\Gamma,\calO)$.
\item
An
\mbox{$\calO$-linear}
map between
$\calO$-lattices
$\varphi\colon (\Gamma, b) \to (\Gamma', b')$
is called {\em orthogonal\/} if
$b'(\varphi(x),\varphi(y)) = b(x,y)$
holds for all
$x,y \in \Gamma$.
And~similarly for
$m\times n$-matrices
if
$\Gamma = \calO^n$
and~$\Gamma' = \calO^m$.
\item
We let
$\bbP$
denote the set of all prime~numbers.
\item
For~$a \in \bbQ_2$,
we let
$\nu_2(a) \in \bbZ \cup \{\infty\}$
be the {\em
\mbox{$2$-adic}
exponential valuation\/}
of~$a$.
I.e.\ the exponent
of~$2$
in the unique factorisation
$a = u2^{\nu_2(a)}$,
for
$u\in\bbZ_2$
a~unit.

Putting
$\nu_2(A) := \min \{\nu_2(a_{ij}) \mid i=1,\ldots,m,\, j = 1,\ldots,n\} \,,$
we extend the
$2$-adic
valuation from
$\bbQ_2$
to matrices
$\smash{A = (a_{ij}) \in \Mat_{m \times n}(\bbQ_2)}$.
\item
A generic line on the projective plane
$\Pb^2$
is denoted
by~$l$.
\end{iii}
\end{notao}

\section{The surfaces studied}

\begin{ttt}
\label{six_lines}
Let
$k$
be a field of characteristic
$\neq \!2$
and
$\smash{l_1, \ldots, l_6 \in \Gamma(\Pb^2_{\overline{k}}, \calO(1))}$
six linear forms, such that the vanishing loci of any three of them do not have a geometric point in~common. Suppose~that
$\{l_1, \ldots, l_6\}$
is a
\mbox{$\smash{\Gal(\overline{k}/k)}$-invariant}
set. Then the double cover
$S'$
of~$\Pb^2_k$,
given by
\begin{equation}
\label{six_lines_eq}
W^2 = l_1 \cdots l_6 \,,
\end{equation}
geometrically has 15 isolated singularities, which are ordinary double points. The~minimal resolution of singularities is a
surface~$S$
of
type~$K3$.
This~is, in fact, one of the most classical families of
$K3$~surfaces,
cf.\ \cite[Chapter 6, Section 2]{GH}.
\end{ttt}

\begin{exs}
\label{histo_exs}
We~actually work
over~$\bbZ$
and consider the double covers
of~$\Pb^2_\bbZ$,
given by the equations~below,

\begin{iii}
\item
$S'_1 \colon W^2 = T_1T_2T_3(T_1+T_2+T_3)(3T_1+5T_2+7T_3)(-5T_1+11T_2-2T_3)$, 
\item
$S'_2 \colon W^2 = T_1T_2T_3(2T_1+4T_2-3T_3)(T_1-5T_2-3T_3)(T_1+3T_2+3T_3)$,  
\item
$S'_3 \colon W^2 = T_1T_2T_3(4T_1+9T_2+T_3)(-T_1-T_2-4T_3)(16T_1+25T_2+T_3)$, 
\item
$S'_4 \colon W^2 = T_1T_2T_3(T_1+T_2+T_3)(T_1+2T_2+3T_3)(5T_1+8T_2+20T_3)$.   
\item
$S'_5 \colon W^2 = T_1T_2 (T_1^4 - 7T_1^3T_2 - T_1^3T_3 + 19T_1^2T_2^2 + 4T_1^2T_2T_3 + T_1^2T_3^2 - 23T_1T_2^3 - 7T_1T_2^2T_3 \hfill$\vspace{-.8mm}

$\hfill{}- 6T_1T_2T_3^2 - T_1T_3^3 + 11T_2^4 + 7T_2^3T_3 + 9T_2^2T_3^2 + 3T_2T_3^3 + T_3^4)$.

Note~that, in the equation
for~$S'_5$,
the homogeneous degree four factor on the right hand side completely splits over the integer ring of the fifth cyclotomic
field~$\smash{\bbQ(\zeta_5)}$,
the~irreducible factors being
$(T_1 - (\zeta_5 + 2)T_2 + \zeta_5^2T_3)$
and its conjugates.
\end{iii}
In~each case, we let
$S_i$
be the blowing-up
of~$S'_i$
in the Zariski closure, equipped with the induced reduced scheme structure, of the singular locus of the generic
fibre~$\smash{S'_{i,\bbQ}}$.
In~particular,
$\smash{S_{i,\bbQ}}$
is the minimal resolution of singularities
of~$\smash{S'_{i,\bbQ}}$.
\end{exs}

\begin{rem}
The surfaces above were investigated in detail in a project concerning the Frobenius trace distributions and the Sato--Tate conjecture for
$K3$~surfaces.
The~efficient point counting algorithm described here was used in order to speed up the computations that led to the histograms presented there~\cite{EJ21}.
\end{rem}

\begin{rems}
\label{main_prop}
Let~us recall from~\cite[Section~5]{EJ21} the main properties of the
surfaces~$S_i$,
for
$i = 1, \ldots, 5$.

\begin{abc}
\item
(Geometric Picard ranks)
In~each case, the pull back of a general line
on~$\smash{\Pb^2_{\overline\bbQ}}$,
together with the exceptional curves resulting from the resolution of singularities, generates a
rank-$16$
sublattice in the geometric Picard~group. The~geometric Picard rank of each of the surfaces
$\smash{S_{1,\bbQ}}$,
$\smash{S_{2,\bbQ}}$,
and~$\smash{S_{4,\bbQ}}$
is indeed equal
to~$16$,
while the
surface~$\smash{S_{3,\bbQ}}$
is of geometric Picard
rank~$17$.
\item
(Bad primes)
Aside from the
prime~$2$,
these are exactly the primes
$p$
such that,
modulo~$p$,
some combination of three of the branch lines has at least one point in common. Thus, the bad primes could be determined by factoring the determinants
$\det (l_i \,\, l_{i'} \,l_{i''})$,
for~$\{i, i'\!, i''\} \subset \{1, \ldots, 6\}$
any subset of size~three. The computation results in the sets
$\{2, 3, 5, 7, 11, 13, 29\}$,
$\{2, 3, 5, 7\}$,
$\{2, 3, 5, 7, 11\}$,
$\{2, 3, 5\}$,
and
$\{2, 5\}$
of bad primes, for
$S_1$,
$S_2$,
$S_3$,
$S_4$
and~$S_5$,
respectively.
\item
(Special properties)
The surface
$S_1$
is a generic surface with a branch locus of six lines, while
$S_2$
has a trivial jump character~\cite{CEJ}. 
The surface
$S_3$
has trivial jump character, too, but higher Picard~rank. Finally,~the surface
$S_4(\bbC)$
has complex multiplication by
$\smash{\bbQ(\sqrt{-1})}$
and
$S_5(\bbC)$
is known to have real multiplication by
$\smash{\bbQ(\sqrt{5})}$.
\end{abc}
\end{rems}

\begin{ttt}[An improvement of Deligne's general bound]
\label{Weil_bound_impr}
As the
$K3$~surfaces
$S_i$,
for
$i=1,2,3,4$,
are of geometric Picard rank
$\geq \!16$
and
$\smash{\Pic(S_{i,\overline\bbQ})}$
is a trivial
$\smash{\Gal(\overline\bbQ/\bbQ)}$-module,
the action
of~$\Frob$
on~$\smash{H^2_\et(S_{\overline\bbF_{\!q}}, \bbZ_2(1))}$
is bound to have the
eigenvalue~$1$
at least of
multiplicity~$16$.
Formula~(\ref{Lefschetz}) therefore actually yields the sharper~estimate
\begin{equation}
\label{DelWeil}
|\#S_i(\bbF_{\!p}) - (p^2+16p+1)| \leq 6p \,.
\end{equation}
Thus,~in each of these cases, it suffices to determine
$(\#S_i(\bbF_{\!p}) \bmod 16p)$.
\end{ttt}

\begin{ttt}
We~explain the
\mbox{$2$-adic}
method in Sections~\ref{EtCohBr}, \ref{overdet}, and~\ref{Algorithm} below, for
$K3$~surfaces
of type~(\ref{six_lines_eq}), of the kind that the linear forms
$l_1, \ldots, l_6$
are defined
over~$\bbZ$.
This~covers Examples~\ref{histo_exs}.i), ii), iii), and~iv), but not~v), which is more~advanced. We~report on the modifications necessary in order to deal with
$S_5$
in Section~\ref{general}.

\end{ttt}

\section{\'Etale cohomology and the Brauer group}
\label{EtCohBr}

\subsubsection*{General\/
$K3$
surfaces}
Let~$S$
be a
$K3$~surface
over a
field~$k$.
Then
$\smash{H^2_\et (S_{\overline{k}}, \bbZ_2(1))}$
is a free
$\bbZ_2$-module
of
rank~$22$.
Since
$\smash{H^3_\et(S_{\overline{k}}, \bbZ_2(1)) = 0}$,
the change of coefficients map
$\smash{H^2_\et(S_{\overline{k}}, \bbZ_2(1)) \!\otimes_{\bbZ_2}\! \bbZ/2^i\bbZ \stackrel{\cong\;}{\to} H^2_\et(S_{\overline{k}}, \mu_{2^i})}$
is an isomorphism, for any
$i\in\bbN$.

Let us recall that there is the natural cup product pairing
\begin{equation}
\label{Poincare}
\langle.,.\rangle\colon H^2_\et(S_{\overline{k}}, \bbZ_2(1)) \times H^2_\et(S_{\overline{k}}, \bbZ_2(1)) \to \bbZ_2
\end{equation}
that is non-degenerate and even perfect, by Poincar\'e duality \cite[Expos\'e~XVIII, Th\'e\-o\-r\`eme 3.2.5]{SGA4}. Moreover,~we have at our disposal the Chern class homomorphism
$\smash{c_1\colon \Pic(S_{\overline{k}}) \to H^2_\et (S_{\overline{k}}, \bbZ_2(1))}$
\cite[Expos\'e~VII, Section~3]{SGA5}, under which the intersection pairing agrees with the cup product pairing.

\begin{defi}
\label{transc_lattice}
Let~$S$
be a
$K3$~surface
over a
field~$k$.

\begin{abc}
\item
Then, we write
$\smash{\P(S_{\overline{k}}, \bbZ_2) := c_1(\Pic(S_{\overline{k}}))}$.
\item
The~orthogonal complement
$\smash{\T(S_{\overline{k}}, \bbZ_2) := \P(S_{\overline{k}}, \bbZ_2)^\perp \subset H^2_\et(S_{\overline{k}}, \bbZ_2(1))}$
is called the {\em transcendental lattice\/}
of~$S$.
Note~that
$\smash{\T(S_{\overline{k}}, \bbZ_2)}$,
as well
as~$\smash{\P(S_{\overline{k}}, \bbZ_2)}$,
is a
$\bbZ_2$-lat\-tice.
\end{abc}
\end{defi}

\noindent
According to this definition,
$\smash{\T(S_{\overline{k}}, \bbZ_2) \subset H^2_\et(S_{\overline{k}}, \bbZ_2(1))}$
clearly has no cotorsion. Let~us note that
$\smash{\P(S_{\overline{k}}, \bbZ_2) \subset H^2_\et(S_{\overline{k}}, \bbZ_2(1))}$
has no cotorsion either. Indeed,~suppose, for a certain
$\calL \in \Pic(S_{\overline{k}})$,
that the Chern class
$\smash{c_1(\calL) \in H^2_\et(S_{\overline{k}}, \bbZ_2(1))}$
is divisible
by~$2$.
Then~$\smash{(c_1 \!\otimes_{\bbZ_2}\! \bbZ/2\bbZ)(\calL) = 0 \in H^2_\et(S_{\overline{k}}, \mu_2)}$
and the exactness of the cohomology sequence
$\smash{\Pic(S_{\overline{k}}) \stackrel{\cdot2\;}{\to} \Pic(S_{\overline{k}}) \to H^2_\et(S_{\overline{k}}, \mu_2)}$,
induced by the Kummer \mbox{sequence,} shows that
$\calL$
is divisible
by~$2$
itself.

Consequently,~one~has
\begin{equation}
\label{perpendicular}
\P(S_{\overline{k}}, \bbZ_2) := \T(S_{\overline{k}}, \bbZ_2)^\perp \,,
\end{equation}
too. Indeed, the cup product pairing~(\ref{Poincare}) is non-degenerate.

\begin{ttt}
Let~$k$
be the finite
field~$\bbF_{\!q}$.
Then,~in terms of the geometric Picard group and the transcendental lattice, the Lefschetz trace formula~(\ref{Lefschetz}) takes the~form
\begin{align}
\label{Lefschetz2}
\#S(\bbF_{\!q}) &= q^2 + \Tr(\Frob\colon \Pic(S_{\overline\bbF_{\!q}}) \!\otimes_\bbZ\! \bbQ \to \Pic(S_{\overline\bbF_{\!q}}) \!\otimes_\bbZ\! \bbQ) q \nonumber \\[-1.2mm]
&\hspace{3.9cm} {}+ \Tr(\Frob\colon \T(S_{\overline\bbF_{\!q}}, \bbZ_2) \to \T(S_{\overline\bbF_{\!q}}, \bbZ_2)) q + 1 \,.
\end{align}
Similarly,~for
$S$
a flat
\mbox{$\bbZ$-scheme}
such that
$S_\bbQ$
is a
$K3$~surface
and a
prime~$p$
of good reduction, one has
\begin{align}
\label{Lefschetz3}
\#S(\bbF_{\!p}) &= p^2 + \Tr(\Frob_p\colon \Pic(S_{\overline\bbQ}) \!\otimes_\bbZ\! \bbQ \to \Pic(S_{\overline\bbQ}) \!\otimes_\bbZ\! \bbQ) p \nonumber \\[-1.2mm]
&\hspace{3.9cm} {}+ \Tr(\Frob_p\colon \T(S_{\overline\bbQ}, \bbZ_2) \to \T(S_{\overline\bbQ}, \bbZ_2)) p + 1 \,.
\end{align}
In~particular, this means
$\#S(\bbF_{\!p}) = p^2 + rp + \Tr(\Frob_p\colon \T(S_{\overline\bbQ}, \bbZ_2) \to \T(S_{\overline\bbQ}, \bbZ_2)) p + 1$
in the case that
$\rk\Pic(S_\bbQ) = \rk\Pic(S_{\overline\bbQ}) = r$.
\end{ttt}

\begin{ttt}
On cohomology with
\mbox{$2$-torsion}
coefficients, the cup product pairing induces a canonical
$\Gal(\overline{k}/k)$-equivariant
isomorphism
\begin{align*}
H^2_\et(S_{\overline{k}}, \mu_2) \stackrel{\cong}{\longrightarrow} &\Hom(H^2_\et(S_{\overline{k}}, \mu_2), \bbZ/2\bbZ) \\
 =\;\, &\Hom(H^2_\et(S_{\overline{k}}, \bbZ_2(1)), \bbZ/2\bbZ) \,.
\end{align*}
Restricting the domain on the right hand side from
$H^2_\et(S_{\overline{k}}, \bbZ_2(1))$
to the transcendental lattice, one obtains a canonical homomorphism
\begin{equation}
\label{Geemen}
H^2_\et(S_{\overline{k}}, \mu_2) \longrightarrow \Hom(\T(S_{\overline{k}}, \bbZ_2), \bbZ/2\bbZ) \,,
\end{equation}
which is surjective, since
$\smash{\T(S_{\overline{k}}, \bbZ_2) \subset H^2_\et(S_{\overline{k}}, \bbZ_2(1))}$
has no cotorsion.
\end{ttt}

\begin{lem}
Let\/~$S$
be a\/
$K3$~surface
over a
field\/~$k$.
Then~the kernel of the homomorphism (\ref{Geemen}) coincides with the image of\/
$\smash{c_1 \!\otimes_{\bbZ_2}\! \bbZ/2\bbZ \colon \Pic(S_{\overline{k}}) \to H^2_\et(S_{\overline{k}}, \mu_2)}$.\smallskip

\noindent
{\bf Proof.}
{\em
``$\supseteq$''
is clear.
``$\subseteq$'':
Let
$\smash{\gamma \in H^2_\et(S_{\overline{k}}, \mu_2)}$
be in the kernel of~(\ref{Geemen}). Lift
$\gamma$
to a class
$\smash{\widetilde\gamma \in H^2_\et(S_{\overline{k}}, \bbZ_2(1))}$.
Then, for every
$\smash{\chi \in \T(S_{\overline{k}}, \bbZ_2)}$,
one has that
$\smash{\langle \widetilde\gamma, \chi\rangle \in \bbZ_2}$
is divisible
by~$2$.
Since
$\smash{\T(S_{\overline{k}}, \bbZ_2)}$
has no cotorsion and the pairing (\ref{Poincare}) is perfect, there exists a class
$\smash{\kappa \in H^2_\et(S_{\overline{k}}, \bbZ_2(1))}$
of the kind that
$\smash{\langle \widetilde\gamma, \chi\rangle = 2 \langle\kappa, \chi\rangle}$,
for any
$\smash{\chi \in \T(S_{\overline{k}}, \bbZ_2)}$.
In other words,
$\smash{\widetilde\gamma - 2\kappa \in \T(S_{\overline{k}}, \bbZ_2)^\perp = \P(S_{\overline{k}}, \bbZ_2)}$,
according to (\ref{perpendicular}), which identifies
$\gamma$
as an element in the image
of~$c_1 \!\otimes_{\bbZ_2}\! \bbZ/2\bbZ$.
}
\eop
\end{lem}

We denote by
$\Br(S_{\overline{k}}) := H^2_\et(S_{\overline{k}}, \bbG_m)$
the geometric Brauer group
of~$S$.
A~standard application of the Kummer sequence shows that the
\mbox{$2$-torsion}
part is given by
$\Br(S_{\overline{k}})_2 = H^2_\et(S_{\overline{k}}, \mu_2) / \im(c_1 \!\otimes_{\bbZ_2}\! \bbZ/2\bbZ)$.

\begin{theo}
\label{prop_Br2}
Let\/~$S$
be a\/
$K3$~surface
over a
field\/~$k$.
Then~the homomorphism~(\ref{Geemen}) induces a canonical\/
$\smash{\Gal(\overline{k}/k)}$-equivariant
isomorphism
\begin{align}
\label{Br2}
\Br(S_{\overline{k}})_2 \stackrel{\cong}{\longrightarrow} \Hom(\T(S_{\overline{k}}, \bbZ_2), \bbZ/2\bbZ) \,.
\end{align}
\end{theo}

\begin{rem}
In~other words, the transcendental lattice
modulo~$2$
is dual to the
\mbox{$2$-torsion}
of the Brauer~group. This~has been known before, at least in the context of complex analytic
$K3$~surfaces,
cf.\ \cite[Paragraph~2.1]{vG}.
\end{rem}

\subsubsection*{Double covers of\/
$\Pb^2$
branched over six lines}
Let~$S$
be a
$K3$~surface
over a
field~$k$
of the kind described in paragraph~\ref{six_lines}.
Then~$S$
may be obtained as a double cover
of~$B$,
the projective plane, blown up in the singular locus of
$V(l_1 \cdots l_6)$,
which forms a reduced
\mbox{$k$-scheme}
of
length~$15$.
The branch locus of the double cover
$\pi\colon S \to B$
is the strict transform
of~$V(l_1 \cdots l_6)$.
This~is a disjoint union of six projective~lines.

\begin{theo}[A.\,N.\ Skorobogatov]
\label{Skoro}
Let\/~$k$
be a field of characteristic
not\/~$2$
and let\/
$S$
be a\/
$K3$~surface\/
over\/~$k$
as in~\ref{six_lines}. Then~there is a\/
$\smash{\Gal(\overline{k}/k)}$-equivariant
isomorphism
$$\Br(S_{\overline{k}})_2 \stackrel{\cong}{\longrightarrow} \Pic(B_{\overline{k}})^\even /\pi_* \Pic(S_{\overline{k}}) \,,$$
for\/
$\Pic(B_{\overline{k}})^\even \subseteq \Pic(B_{\overline{k}})$
the subgroup formed by the classes having an even intersection number with each connected component of the branch~locus.\smallskip

\noindent
{\bf Proof.}
{\em
This~is a particular case of A.\,N.~Skorobogatov's explicit description of
$\Br(S_{\overline{k}})_2$
for double covers~\cite[Theorem~1.1]{Sk}. Note~that the geometric Picard group of the branch locus has no~torsion.
}
\eop
\end{theo}

\subsubsection*{The case of six\/
\mbox{$k$-rational}
lines}

\begin{coro}
\label{Gal_triv_T2T}
Let\/~$k$
be a field of characteristic
not\/~$2$
and\/
$S$
a\/
$K3$~surface
over\/~$k$
as in~\ref{six_lines}. Suppose~that\/
$l_1, \ldots, l_6$
are defined
over\/~$k$.

\begin{abc}
\item
Then~the natural\/
$\smash{\Gal(\overline{k}/k)}$-action
on\/~$\smash{\Br(S_{\overline{k}})_2}$
is~trivial.
\item
The~natural\/
$\smash{\Gal(\overline{k}/k)}$-action
on\/~$\smash{\T(S_{\overline{k}}, \bbZ_2)/2\T(S_{\overline{k}}, \bbZ_2)}$
is trivial,~too.
\end{abc}\smallskip

\noindent\looseness-1
{\bf Proof.}
{\em
a)
Since~$l_1, \ldots, l_6$
are defined
over~$k$,
$B$
is the blowing-up of
$\smash{\Pb^2_k}$
in
$15$
\mbox{$k$-rational}
points. Therefore,~the
$\smash{\Gal(\overline{k}/k)}$-action
on the whole of
$\smash{\Pic(B_{\overline{k}})}$
is~trivial.\smallskip

\noindent
b)
follows from~a), together with Theorem~\ref{prop_Br2}.
}
\eop
\end{coro}

\begin{coro}[The splitting field of
$\T(S_{\overline{k}}, \bbZ_2)/4\T(S_{\overline{k}}, \bbZ_2)$]
\label{split_T4}
Let\/~$k$
be a field of characteristic
not\/~$2$
and\/
$S$
a\/
$K3$~surface
over\/~$k$
as in~\ref{six_lines}. Suppose~that\/
$l_1, \ldots, l_6$
are defined
over\/~$k$.
Denote
by\/~$K \supseteq k$
the splitting field of\/
$\T(S_{\overline{k}}, \bbZ_2)/4\T(S_{\overline{k}}, \bbZ_2)$.

\begin{abc}
\item
Then\/~$K$
is an abelian extension
of\/~$k$
of exponent at
most\/~$2$.
\item
Suppose~that\/
$k$
is a number field.
Then\/~$K$
is unramified
over\/~$k$
at all primes of good reduction and odd residue~characteristic.
\end{abc}\smallskip

\noindent
{\bf Proof.}
{\em
a)
By~definition, one has a natural injection
$$\Gal(K/k) \hookrightarrow \Aut(\T(S_{\overline{k}}, \bbZ_2)/4\T(S_{\overline{k}}, \bbZ_2)) \,.$$
But~$\{A \in \GL_n(\bbZ / 4 \bbZ) \mid A \equiv E_n \pmodulo 2\}$
is an elementary abelian
\mbox{$2$-group},
for any
$n\in\bbN$.\smallskip

\noindent
b)
As
$\T(S_{\overline{k}}, \bbZ_2) \subset H^2_\et(S_{\overline{k}}, \bbZ_2(1))$
has no cotorsion, the natural homomorphism
$\T(S_{\overline{k}}, \bbZ_2)/4\T(S_{\overline{k}}, \bbZ_2) \hookrightarrow H^2_\et(S_{\overline{k}}, \mu_4)$
is~injective. Moreover,~by virtue of the smooth specialisation theorem~\cite[Expos\'e XVI, Corollaire~2.3]{SGA4}, the splitting field of
$H^2_\et(S_{\overline{k}}, \mu_4)$
is known to be unramified at any prime
of~$k$
of odd residue characteristic, at which
$S$~has
good~reduction.
}
\eop
\end{coro}

\begin{coro}
\label{Fall_Q}
Let\/~$S$
be a\/
$K3$~surface
as in~\ref{six_lines},
over\/~$k=\bbQ$.
Suppose~that\/
$l_1, \ldots, l_6$
are defined
over\/~$\bbQ$.
Then,~for an odd
prime\/~$p$
of good reduction, the action of the Frobenius\/
$\Frob_p$
on\/
$\T(S_{\overline\bbQ}, \bbZ_2)/4\T(S_{\overline\bbQ}, \bbZ_2)$
is completely determined by the class of\/
$\Frob_p$
in the Galois group of the number~field
$$K = \bbQ(\sqrt{-1}, \sqrt{2}, \sqrt{d} \mid d \mbox{ \rm a bad prime for } S) \,.$$
{\bf Proof.}
{\em
This~is the particular case of Corollary~\ref{split_T4},
for~$k = \bbQ$.
}
\eop
\end{coro}

\subsubsection*{The action
of\/~$\Sym(6)$
on\/~$\T(S_{\overline{k}}, \bbZ_2)/2\T(S_{\overline{k}}, \bbZ_2)$
in the case that\/
$\rk\Pic(S_{\overline{k}}) = 16$}

\begin{nota}
In~the remainder of this section, a double index is meant to be an unordered pair. E.g.,
$a_{ij} = a_{ji}$.
In~particular, we write
$e_{ij}$,
but also
$e_{ji}$,
for the exceptional curve
on~$B_{\overline{k}}$
that lies over the point of intersection
$\V(l_i) \cap \V(l_j)$,
for~$i \neq j$,
$1 \leq i, j \leq 6$.
\end{nota}

\begin{lem}
\label{Pic_exakt}
Let\/~$k$
be a field of characteristic
not\/~$2$
and\/
$S$
a\/
$K3$~surface
over\/~$k$
as in~\ref{six_lines}. Suppose~that\/
$\rk\Pic(S_{\overline{k}}) = 16$.

\begin{abc}
\item[{\rm a.i) }]
Then\/
$\Pic(S_{\overline{k}}) \supset \bbZ\pi^*[l] \oplus \bbZ\pi^*[e_{12}] \oplus\cdots\oplus \bbZ\pi^*[e_{56}]$
is a sublattice of full~rank.
\item[{\rm ii) }]
There~are further divisor classes\/
$[D_1], \ldots, [D_6] \in \Pic(S_{\overline{k}})$
such~that
\begin{equation}
\label{halbierbar}
2[D_i] = \pi^*[l] + \!\!\sum_{\atop{j=1,\ldots,6,}{j \neq i}} \!\!\!\pi^*[e_{ij}] \,.
\end{equation}
The~classes\/
$[D_1], \ldots, [D_6]$
generate\/~$\Pic(S_{\overline{k}})$,
together with\/
$\pi^*[l], \pi^*[e_{12}], \ldots, \pi^*[e_{56}]$.
\item[{\rm b) }]
Put
\begin{align*}
M := \{a \pi^*[l] + a_{12} \pi^*[e_{12}] + \cdots + a_{56} \pi^*[e_{56}] \in \bbF_{\!2}\pi^*[l] \oplus \bbF_{\!2}\pi^*[e_{12}] \oplus\cdots\oplus \bbF_{\!2}\pi^*[e_{56}]& \\[-1mm]
\textstyle \mid a + \!\!\!\sum\limits_{\atop{j=1,\ldots,6,}{j \neq i}}\!\!\!\! a_{ij} = 0, \text{\rm ~for~} i = 1, \ldots, 6&\} \,.
\end{align*}\vskip-3.2mm
Then~there is a natural isomorphism
\begin{equation}
\label{Br2_expl}
\textstyle \Br(S_{\overline{k}})_2 \cong M/\langle\pi^*[l] + \!\!\!\sum\limits_{\atop{j=1,\ldots,6,}{j \neq i}} \!\!\!\!\pi^*[e_{ij}] \mid i = 1, \ldots, 6\rangle \,.
\end{equation}
\end{abc}

\noindent
{\bf Proof.}
{\em
One has
$\Pic(B_{\overline{k}}) = \bbZ[l] \oplus \bbZ[e_{12}] \oplus\cdots\oplus \bbZ[e_{56}]$,
the direct sum being orthogonal. Moreover,
$[l] \!\cdot\! [l] = 1$
and
$[e_{ij}] \!\cdot\! [e_{ij}] = -1$,
for
$1 \leq i < j \leq 6$.
As
$\pi_{\overline{k}}\colon S_{\overline{k}} \to B_{\overline{k}}$
is finite of
degree~$2$,
this yields that
$\pi^*[l], \pi^*[e_{12}], \ldots, \pi^*[e_{56}] \in \Pic(S_{\overline{k}})$
are mutually perpendicular, with
$\pi^*[l] \!\cdot\! \pi^*[l] = 2$
and~$\pi^*[e_{ij}] \!\cdot\! \pi^*[e_{ij}] = -2$.
From~this, a.i) immediately~follows.
Consequently,~one has
$\bbQ\pi^*[l] \oplus \bbQ\pi^*[e_{12}] \oplus\cdots\oplus \bbQ\pi^*[e_{56}] \supset \Pic(S_{\overline{k}})$.
As~the intersection numbers with the base elements have to be integers, the coefficients are in fact half integral, at~most.

Furthermore,~over the strict transform of the quintic
$\smash{\V(l_2 \cdots l_6 - l_1^5) \subset \Pb^2_{\overline{k}}}$,
the double cover
$\smash{\pi_{\overline{k}}\colon S_{\overline{k}} \to B_{\overline{k}}}$
splits, since the equation of the surface goes over into
$w^2 = l_1^6$
\cite[Remark~4.6]{EJ18}. This~yields a divisor
$D \in \Div(S_{\overline{k}})$
such that
$\pi_*[D] = 5[l] + [e_{12}] + [e_{13}] + [e_{14}] + [e_{15}] + [e_{16}]$.
One~may put
$[D_1] := [D] - 2\pi^*[l]$
in order to fulfil~(\ref{halbierbar}). The divisor classes
$[D_2], \ldots, [D_6]$
are constructed~analogously.
Thus,~there is an~inclusion
\begin{equation}
\label{Pic_incl}
\Pic(B_{\overline{k}}) \supseteq (\bbZ[l] \oplus \bbZ[e_{12}] \oplus\cdots\oplus \bbZ[e_{56}]) + [D_1] + \cdots + [D_6] \,,
\end{equation}
and in order to complete the proof of~a), it needs to be shown that equality~holds.

On~the other hand, one has
$\smash{M = \Pic(B_{\overline{k}})^\even/\pi_* (\bbZ\pi^*[l] \!\oplus\! \bbZ\pi^*[e_{12}] \!\oplus\!\cdots\!\oplus\! \bbZ\pi^*[e_{56}])}$.
Hence,~by Theorem~\ref{Skoro}, there is actually a natural surjection
\begin{equation}
\label{Br2_surj}
\textstyle \Br(S_{\overline{k}})_2 \twoheadrightarrow M/\langle\pi^*[l] + \!\!\sum\limits_{\atop{j=1,\ldots,6,}{j \neq i}} \!\!\!\pi^*[e_{ij}] \mid i = 1, \ldots, 6\rangle \vspace{-1.5mm}
\end{equation}
that is a bijection if and only if equality holds in~(\ref{Pic_incl}). But~an explicit calculation reveals that the right hand side of~(\ref{Br2_surj}) is of
order~$64$,
$M$
being of
order~$2048$.
As,~by Theorem~\ref{prop_Br2},
$\smash{\Br(S_{\overline{k}})_2}$
is of
order~$64$,
too, this proves both, a) and~b).%
}%
\eop
\end{lem}

\begin{theo}
\label{S6_action}
Let\/~$k$
be a field of characteristic
not\/~$2$
and\/
$S$
a\/
$K3$~surface
over\/~$k$
as in~\ref{six_lines}. Suppose~that\/
$\smash{\rk\Pic(S_{\overline{k}}) = 16}$.

\begin{abc}
\item[{\rm a.i) }]
Then~the group\/
$\Sym(6)$
permuting the six branch lines naturally acts on the sublattice of\/
$\smash{\Pic(S_{\overline{k}})}$
described in Lemma~\ref{Pic_exakt}.a.i)~by
$$\sigma(\pi^*[l]) = \pi^*[l] \qquad\text{and}\qquad \sigma(\pi^*[e_{ij}]) = \pi^*[e_{\sigma(i)\sigma(j)}] \,.$$
\item[{\rm ii) }]
The action on\/
$\smash{\Br(S_{\overline{k}})_2}$
is as~follows. There~is an\/
\mbox{$\bbF_{\!2}$-basis\/}
$(b_1, \ldots, b_6)$
of\/~$\smash{\Br(S_{\overline{k}})_2}$
such that
$$\sigma(b_i) = \left\{
\begin{array}{lc}
b_{\widetilde\sigma(i)} \,,     & {\rm ~if\;} \sigma\in\Alt(6) \,, \\
b_{\widetilde\sigma(i)} + c \,, & {\rm ~if\;} \sigma\not\in\Alt(6) \,,
\end{array}
\right.$$
for\/~$c := b_1 + \cdots + b_6$.
Here,
$\Sym(6) \to \Sym(6)$,
$\sigma \mapsto \widetilde\sigma$,
is an outer automorphism.
\item[{\rm b) }]
The~natural\/
\mbox{$\Gal(\overline{k}/k)$-actions}
on\/~$\smash{\Pic(S_{\overline{k}})}$
and\/~$\smash{\Br(S_{\overline{k}})_2}$
are the compositions of the natural\/
\mbox{$\Gal(\overline{k}/k)$-action}
on the six branch lines with the actions described in~a).
\end{abc}\smallskip

\noindent
{\bf Proof.}
{\em
a.i) and~b) are clear.\smallskip

\noindent
a.ii)
One~puts
$\smash{b_6 := \overline{\pi^*[e_{12}] + \pi^*[e_{23}] + \pi^*[e_{34}] + \pi^*[e_{45}] + \pi^*[e_{15}]} \in \Br(S_{\overline{k}})_2}$.~Note
that indeed
$\pi^*[e_{12}] + \pi^*[e_{23}] + \pi^*[e_{34}] + \pi^*[e_{45}] + \pi^*[e_{15}] \in M$,
so that this is a correct definition.

It~is clear that
$b_6$
is stabilised by a dihedral group of order ten, permuting only
$\{1, \ldots, 5\}$.
Moreover,~applying the relations in (\ref{Br2_expl}), for
$i=1$
and~$4$,
one finds that
$\smash{b_6 = \overline{\pi^*[e_{13}] + \pi^*[e_{23}] + \pi^*[e_{24}] + \pi^*[e_{46}] + \pi^*[e_{16}]}}$,
too. Hence,~the stabiliser
of~$b_6$
is a
\mbox{$2$-transitive}
subgroup
of~$\Sym(6)$,
of order a multiple
of~$60$.

A~machine calculation shows that the orbit
of~$b_6$
under~$\Sym(6)$
is indeed of size~twelve, so that the stabiliser
of~$b_6$
is
$\Alt(5)$,
transitively embedded
into~$\Sym(6)$
\cite[Table~2.1]{DM}. Furthermore,~the orbit
$\{b_1, \ldots, b_6\}$
of~$b_6$
under~$\Alt(6) \subset \Sym(6)$
turns out to be
\mbox{$\bbF_{\!2}$-linearly}
independent.

The~orbit
under~$\Sym(6)$
is, in fact,
$\{b_1, \ldots, b_6, b_1+c, \ldots, b_6+c\}$.
Moreover,
$\smash{c \in \Br(S_{\overline{k}})_2 \!\setminus\! \{0\}}$
is the unique
\mbox{$\Alt(6)$-invariant}
element, and hence
$\Sym(6)$-invariant.
Consequently,~the stabiliser of the
class~$\smash{\overline{b}_6 \in \Br(S_{\overline{k}})_2/\langle c\rangle}$
is
\mbox{$2$-transitive}
of
order~$120$,
and thus
$\Sym(5)$,
transitively embedded
into~$\Sym(6)$.
This~shows that indeed the~classes
$\smash{\overline{b}_1, \ldots, \overline{b}_6 \in \Br(S_{\overline{k}})_2/\langle c\rangle}$
are permuted according to an outer automorphism of the group
$\Sym(6)$
\cite[Subsection~2.4.2]{Wi}, which completes the~proof.
}
\eop
\end{theo}

\section{A theorem on the
$2$-adic
orthogonal group}
\label{overdet}

\subsubsection*{The main theorem.
$2$-adic
overdetermination of the trace}

\begin{theo}
\label{Overdet}
Let\/~$n\in\bbN$.
With~respect to a non-degenerate, symmetric bilinear form
on\/~$\bbQ_2^n$,
let\/
$U_1, U_2 \in \Mat_{n \times n}(\bbZ_2)$
be orthogonal matrices such~that
$$U_1 \equiv U_2 \pmodulo 4 \,.$$

\begin{abc}
\item
If\/
$U_1 \equiv E_n \pmodulo 2$
then\/~$\Tr(U_1) \equiv \Tr(U_2) \pmodulo {16}$.
\item
If\/
$U_1^2 \equiv E_n \pmodulo 2$
then\/~$\Tr(U_1) \equiv \Tr(U_2) \pmodulo 8$.
\end{abc}\smallskip

\noindent
{\bf Proof.}
{\em
One~has
$\det(U_1) = \pm1$,
since
$U_1$
is orthogonal. Hence,
$U_1^{-1} \in \Mat_{n \times n}(\bbZ_2)$.
The~congruence
$U_1 \equiv U_2 \pmodulo 4$
therefore implies that
$U_1^{-1} U_2 \equiv E_n \pmodulo 4$.
I.e.,~there exists a matrix
$B \in \Mat_{n \times n}(\bbZ_2)$
such~that
$$U_2 = U_1 (E_n + 4B) \,.$$
Moreover,
$\det(E_n + 4B) = \det(U_1^{-1} U_2) = 1$,
hence the linear approximation of
$\det$
near the unit matrix~yields
$1 = \det(E_n + 4B) \equiv \det(E_n) + 4\Tr(B) \pmodulo{16}$.
I.e.,
\begin{equation}
\label{B4}
\Tr(B) \equiv 0 \pmodulo 4 \,.
\end{equation}

\noindent
a)
Writing
$U_1 = E_n + 2A$,
one~finds
\begin{align*}
\Tr(U_2) = \Tr((E_n + 2A) \!\cdot\! (E_n + 4B)) &= \Tr(E_n + 2A) + 4\Tr(B) + 8 \Tr(AB) \\
 &\equiv \Tr(U_1) + 8 \Tr(AB) \pmodulo{16} \,.
\end{align*}
Finally,~one has
$\Tr(AB) \equiv 0 \pmodulo 2$,
due to Theorem~\ref{Hauptsatz}.a),~below.\smallskip

\noindent
b)
Analogously,~writing
$U_1 = E_n + A$,
one~sees
$$\Tr(U_2) = \Tr((E_n + A) \!\cdot\! (E_n + 4B)) = \Tr(U_1) + 4\Tr(B) + 4\Tr(AB) \,,$$
so that the assertion follows from (\ref{B4}) and Theorem~\ref{Hauptsatz}.b).
}
\eop
\end{theo}

\subsubsection*{$2$-adic
divisibility of traces.}

\begin{theo}
\label{Hauptsatz}
For\/~$n\in\bbN$,
let\/
$b$
be a non-degenerate, symmetric bilinear form
on\/~$\bbQ_2^n$.

\begin{abc}
\item
Let\/
$A,B \in \Mat_{n \times n}(\bbZ_2)$
be such that\/
$E_n + 2A$
and\/~$E_n + 4B$
are orthogonal with respect
to\/~$b$.
Then\/
$\Tr(AB) \equiv 0 \pmod 2$.
\item
Let\/
$A,B \in \Mat_{n \times n}(\bbZ_2)$
be such that\/
$E_n + A$
and\/
$E_n + 4B$
are orthogonal with respect
to\/~$b$
and\/
$(E_n + A)^2 \equiv E_n \pmod 2$.
Then\/
$\Tr(A B) \equiv 0 \pmod 2$.
\end{abc}
\end{theo}

A more natural formulation, not explicitly considering matrices, goes as~follows.

\begin{theo}
\label{trace_div}
Let\/~$\Gamma$
be a\/
\mbox{$\bbZ_2$-lattice}
and\/
$\varphi,\psi\colon \Gamma \to \Gamma$
orthogonal~endomorphisms.

\begin{abc}
\item
If\/
$(\varphi \bmod 2)$
acts trivially on\/
$\Gamma/2\Gamma$
and\/
$(\psi \bmod 4)$
acts trivially on\/
$\Gamma/4\Gamma$
then\/
$\Tr((\varphi-\id) \!\circ\! (\psi-\id)) \equiv 0 \pmodulo {16}$.
\item
If\/
$(\varphi\!\circ\!\varphi \bmod 2)$
acts trivially on\/
$\Gamma/2\Gamma$
and\/
$(\psi \bmod 4)$
acts trivially on\/
$\Gamma/4\Gamma$
then\/~$\Tr((\varphi-\id) \!\circ\! (\psi-\id)) \equiv 0 \pmodulo 8$.
\end{abc}
\end{theo}

\begin{ttt}[Structure of the proof]
The proof of Theorem~\ref{trace_div} is organised as follows.

\begin{iii}
\item
We~directly check the analogue of Theorem~\ref{trace_div} for regular
$\bbZ_2[\sqrt{2}]$-lattices.
\item
We~show Theorem~\ref{trace_div} in a particular case, in which the discriminant of
$b$
has a small
\mbox{$2$-adic}
valuation. The~idea is to construct a
$\smash{\bbZ_2[\sqrt{2}]}$-lattice
$\smash{\widetilde\Gamma \supseteq \Gamma \!\otimes_{\bbZ_2}\! \bbZ_2[\sqrt{2}]}$
such that
$b$
is regular on
$\smash{\widetilde\Gamma}$.
We~check that all relevant properties of the endomorphisms
$\varphi$
and~$\psi$
are preserved under this change of~lattices.
\item
For~the general case, the idea is as follows. By inspection of the dual lattice
$\Gamma^\vee$,
we show that the endomorphisms
$\varphi$
and~$\psi$
respect various other lattices, as~well. We~choose one rather particular such
lattice~$\smash{\Gamma' \subset \Gamma \!\otimes_{\bbZ_2}\! \bbQ_2}$,
thereby making sure that the discriminant of
$b$
on
$\Gamma'$
has a
\mbox{$2$-adic}
valuation sufficiently small, so that ii) applies.
\end{iii}
\end{ttt}

\begin{rems}
\begin{abc}
\item
For the concept of an
$\calO$-lattice,
in general, recall the terminology and notation fixed in the introduction.
\item
The ring
$\smash{\calO = \bbZ_2[\sqrt{2}]}$
is a principal ideal domain, in fact a discrete valuation~ring. In~particular, every finitely generated
$\smash{\calO = \bbZ_2[\sqrt{2}]}$-module
that is projective or only torsion-free is automatically~free.
\item
There~is a more elementary proof for Theorem~\ref{trace_div} that avoids
$\smash{\bbZ_2[\sqrt{2}]}$-lattices
at the cost of a more complicated case distinction.
\end{abc}
\end{rems}

\subsubsection*{First step of the proof -- A variant for regular lattices over discrete valuation rings}

\begin{prop}
\label{Z2_var}
Let\/~$\calO$
be a discrete valuation ring, in which\/
$2 \neq 0$
is not a unit,
$\Gamma$
a regular\/
$\smash{\calO}$-lattice,
and\/
$\varphi,\psi\colon \Gamma \to \Gamma$
orthogonal~endomorphisms.

\begin{abc}
\item
If\/
$(\varphi \bmod 2)$
acts trivially on\/
$\Gamma/2\Gamma$
and\/
$(\psi \bmod 4)$
acts trivially on\/
$\Gamma/4\Gamma$
then\/
$\Tr((\varphi-\id) \!\circ\! (\psi-\id)) \equiv 0 \pmodulo {16}$.
\item
If\/
$(\varphi\circ\varphi \bmod 2)$
acts trivially on\/
$\Gamma/2\Gamma$
and\/
$(\psi \bmod 4)$
acts trivially on\/
$\Gamma/4\Gamma$
then\/
$\Tr((\varphi-\id) \!\circ\! (\psi-\id)) \equiv 0 \pmodulo 8$.
\end{abc}\smallskip

\noindent
{\bf Proof.}
{\em
Choose a basis
for~$\Gamma$
and denote the rank
of~$\Gamma$
by~$n$.
The~bilinear form
on~$\Gamma$
is then given by a symmetric
$n \times n$
matrix~$M$.
Thanks~to the regularity assumption, one has
$\smash{M \in \GL_n(\calO)}$.
Moreover,~the matrix representing
$\psi$
can be written as
$E_n + 4B$,
for some
$\smash{B \in \Mat_{n\times n}(\calO)}$.
By~Lemma~\ref{Lemma_antisym}.b),
$MB$
is symmetric
modulo~$2$
all diagonal coefficients being divisible
by~$2$.\smallskip

\noindent
a)
Here, the matrix representing
$\varphi$
is
$E_n + 2A$,
for a certain
$\smash{A \in \Mat_{n\times n}(\calO)}$.
Thus,~the assertion follows from Lemmas~\ref{Lemma_antisym}.a) and~\ref{Letztes_Lemma}, below.\smallskip

\noindent
b)
Let
$\smash{A \in \Mat_{n\times n}(\calO)}$
be the matrix representing
$\varphi - \id$.
The~fact that
$\varphi$
is orthogonal is then equivalent to
\begin{equation}
\label{orthA}
A^\top M + M A + A^\top M A = 0 \,,
\end{equation}
which implies
$$A^\top M A + M A^2 + A^\top M A^2 = 0 \,.$$
On~the other hand, the assumption
$(E_n+A)^2 \equiv E_n \pmodulo 2$
yields
$A^2 \equiv 0 \pmodulo 2$,
so that one has
$A^\top M A \equiv 0 \pmodulo 2$.
Therefore, formula~(\ref{orthA}), together with the fact that
$M$~is
symmetric, yields that
$MA$
is symmetric
modulo~$2$.
To~summarise, it is proved that Lemma~\ref{Letztes_Lemma} is applicable, which shows that
$\Tr(AB)$
is divisible
by~2.
}
\eop
\end{prop}

\begin{lem}
\label{Lemma_antisym}
Let\/~$\calO$
be a discrete valuation ring, in which\/
$2 \neq 0$
is not a unit, and\/
$M \in \GL_n(\calO)$
a symmetric~matrix.

\begin{abc}
\item
Let\/
$C \in \Mat_{n \times n}(\calO)$
be such that\/
$E_n + 2C$
is orthogonal with respect
to\/~$M$.
Then~the reduction
modulo\/~$2$
of\/
$MC$
is~symmetric.
\item
Let\/
$C \in \Mat_{n \times n}(\calO)$
be such that\/
$E_n + 4C$
is orthogonal with respect
to\/~$M$.
Then~the reduction
modulo\/~$2$
of\/
$MC$
is~symmetric. Furthermore, every diagonal coefficient of\/
$MC$
is divisible
by\/~$2$.
\end{abc}\smallskip

\noindent
{\bf Proof.}
{\em
a)
The orthogonality condition explicitly~reads
\begin{equation}
\label{orth}
C^\top M + M C + 2 C^\top M C = 0 \,.
\end{equation}
Moreover,~as
$M$
is symmetric, we have
$C^\top M = (MC)^\top$.
Thus, (\ref{orth}) implies that
$MC \equiv (MC)^\top \pmodulo 2$,
which shows~a).\smallskip

\noindent
b)
Here, orthogonality means
$C^\top M + M C + 4 C^\top M C = 0$.
Moreover,~once again, we have
$C^\top M = (MC)^\top$,
so that
$MC \equiv -(MC)^\top \pmodulo 4$
follows. This~proves both conclusions of~b).
}
\eop
\end{lem}

\begin{lem}
\label{Letztes_Lemma}
Let\/~$\calO$
be a discrete valuation ring, in which\/
$2 \neq 0$
is not a unit,
$M \in \GL_n(\calO)$
symmetric, and\/
$A, B \in \Mat_{n \times n}(\calO)$
such that\/
$\smash{MA \equiv (MA)^\top \pmodulo 2}$,
$MB \equiv (MB)^\top \pmodulo 2$,
and all diagonal coefficients of
$MB$\/
are divisible\/
by~$2$.
Then\/
$\Tr(AB)$
is divisible
by\/~$2$.\smallskip

\noindent
{\bf Proof.}
{\em
Writing
$U := A M^{-1}$
and
$V := MB$,
one has
\begin{equation}
\label{Trick}
\Tr(AB) = \Tr(A M^{-1} M B) = \Tr(UV) = \sum_{i,j} u_{ij} v_{ji} \,. \vspace{-1.5mm}
\end{equation}
Here,
$U = M^{-1} MA M^{-1}$
is symmetric
modulo~$2$,
since both,
$MA$
and~$M^{-1}$,
are. Furthermore,
$V = MB$
is symmetric
modulo~$2$,
by~assumption.

Therefore,~in (\ref{Trick}), the summands for the indices
$(i,j)$
and
$(j,i)$
coincide
modulo~$2$,
so that the sum of each pair is divisible
by~$2$.
Finally,~the summands for
$i=j$
are divisible
by~$2$,
as the diagonal coefficients
$v_{ii}$
of
$V$~are.
}
\eop
\end{lem}

\subsubsection*{Second step of the proof -- Generalities on
$\bbZ_2$-lattices}

\begin{prop}[Decomposition of
\mbox{$\bbZ_2$-lattices}]
\label{decomp}
Let\/
$(\Gamma, b)$
be a\/
\mbox{$\bbZ_2$-lattice}.
Then there is a decomposition\/
$\smash{\Gamma = \bigoplus_{i = 0}^N \Gamma_i}$
into an orthogonal direct sum of the kind~that
\begin{equation}
\label{Kneser}
b = \bigoplus_{i=0}^N 2^i b_i \,,
\end{equation}
for regular symmetric bilinear forms\/
$b_0, \ldots, b_N$
on\/
$\Gamma_0, \ldots, \Gamma_N$,
respectively.\smallskip

\noindent
{\bf Proof.}
{\em
See~\cite[Chapter~15, Theorem~2]{CS} or~\cite[Satz~15.1]{Kn}. It~is, in fact, shown that, for each
$i$,
the lattice
$\Gamma_i$
may be decomposed further into an orthogonal direct sum of only
$1$-
and
$2$-dimensional
lattices.
}
\eop
\end{prop}

\begin{defi}
Let
$(\Gamma, b)$
be a
\mbox{$\bbZ_2$-lattice}.
Then, by the {\em dual lattice,} we mean
$$\Gamma^\vee := \{x \in \Gamma \!\otimes_{\bbZ_2}\! \bbQ_2 \mid b(x, \Gamma) \subseteq \bbZ_2 \} \,.$$
Here,~the
\mbox{$\bbQ_2$-bilinear}
extension of
$b\colon \Gamma \times \Gamma \to \bbZ_2$
to
$\Gamma \!\otimes_{\bbZ_2}\! \bbQ_2$
is again denoted
by~$b$.
\end{defi}

\begin{ex}
A~decomposition
$\smash{\Gamma = \bigoplus\limits_{i = 0}^N \Gamma_i}$
as above, with
$\smash{b = \bigoplus\limits_{i=0}^N 2^i b_i}$,
yields\smallskip
$$\Gamma^\vee = \bigoplus_{i = 0}^N 2^{-i} \Gamma_i \,.$$
Indeed,~for
each~$i$,
the lattice
$(\Gamma_i, b_i)$
is~regular.
\end{ex}

\begin{lem}
\label{dual}
Let\/~$\Gamma$
be a\/
\mbox{$\bbZ_2$-lattice}
and\/
$\varphi\colon \Gamma \to \Gamma$
an orthogonal~map. Denote~the\/
\mbox{$\bbQ_2$-linear}
extension
of\/~$\varphi$
to\/
$\Gamma \!\otimes_{\bbZ_2}\! \bbQ_2$
again
by\/~$\varphi$.
Then\/
$\varphi(\Gamma^\vee) = \Gamma^\vee$.\smallskip

\noindent
{\bf Proof.}
{\em
We have
\begin{equation}
x \in \Gamma^\vee \Leftrightarrow b(x, \Gamma) \subseteq \bbZ_2 \Leftrightarrow
b(\varphi(x), \varphi(\Gamma)) \subseteq \bbZ_2
\Leftrightarrow
b(\varphi(x), \Gamma) \subseteq \bbZ_2
\Leftrightarrow
\varphi(x) \in \Gamma^\vee \,.
\tag*{$\square$}\vspace{1mm}
\end{equation}
}
\end{lem}

\begin{lem}
\label{Block_Bewertung}
Let\/~$\Gamma$
be a\/
\mbox{$\bbZ_2$-lattice}
and\/
$\varphi\colon\Gamma \to \Gamma$
an orthogonal endomorphism. Moreover,~let\/
$\smash{\Gamma = \bigoplus_{i = 0}^N \Gamma_i}$,
with\/
$\smash{b = \bigoplus_{i=0}^N 2^i b_i}$,
be a decomposition as above and let\/
$\calB$~be
a basis
of\/
$\Gamma$
obtained by concatenating bases of\/
$\Gamma_0, \ldots, \Gamma_N$.
Finally,~put
$$
D = (D_{ij})_{i,j = 0,\ldots,N} := \Mat^\calB_\calB(\varphi) \,,
$$
so that, for\/
$i,j = 0, \ldots, N$,
the block\/
$D_{ij}$
represents an element of\/
$\Hom(\Gamma_j,\Gamma_i)$.

\begin{abc}
\item
Then,~for\/
$i,j = 0, \ldots, N$,
one~has\/
$\nu_2(D_{ij}) \geq j-i$.
\item
Let\/~$e\in\bbN$
and suppose, in addition, that
$\varphi$
acts trivially on the quotient\/
$\Gamma / 2^e \Gamma$.
Then\/
$\nu_2(D_{ij}) \geq j - i + e$,
for\/~$i \neq j$,
$i,j = 0, \ldots, N$.
\end{abc}\smallskip

\noindent
{\bf Proof.}
{\em
a)
As
$\varphi$
maps
$\Gamma$
onto itself, by Lemma~\ref{dual}, it does the same to the dual lattice
$\smash{\Gamma^\vee = \bigoplus_{i = 0}^N 2^{-i} \Gamma_i}$.
This shows that
$\nu_2(D_{ij}) \geq j - i$,
for~$j > i$.
Note that the assertion is trivial in the case
that~$j \leq i$.\smallskip

\noindent
b)
By assumption,
$\varphi$
maps the lattices
$\Gamma$
and
$2^e \Gamma$
onto~themselves. Furthermore,~the action induced on the quotient
$\Gamma / 2^e \Gamma$
is assumed to be~trivial. Consequently,
$\varphi$
maps every lattice
$\Delta$
of the kind that
$2^e \Gamma \subseteq \Delta \subseteq \Gamma$
onto~itself.

We apply this observation to the lattices, given~by
$$\Delta_j := \Bigg(\!\bigoplus_{\atop{i=0, \ldots, N}{i \neq j}} \!\!\!\Gamma_i\!\Bigg)
 \oplus 2^e \Gamma_j \,,$$
for~$j = 0, \ldots, N$.
As~$\varphi$
maps
$\Delta_j$
onto itself, the same is true
for~$\Delta_j^\vee$.
Noticing~that
$$\Delta_j^\vee = \Bigg(\!\!\bigoplus_{\atop{i=0, \ldots, N}{i \neq j}} \!\!\!2^{-i} \Gamma_i\!\Bigg)
 \oplus 2^{-j-e} \Gamma_j \,,$$
this yields
$\nu_2(D_{ij}) \geq j - i + e$,
for~$j \neq i$.
Indeed, one has
that~$\varphi(2^{-j-e} \Gamma_j) \subseteq \Delta_j^\vee$.%
}%
\eop
\end{lem}

\begin{rem}
An alternative proof for~b) may be given as~follows. As
$\varphi$
acts trivially on
$\Gamma / 2^e \Gamma$,
the adjoint map acts trivially on
$\Gamma^\vee / 2^e \Gamma^\vee$.
As the adjoint coincides with the inverse
of~$\varphi$,
the map being orthogonal, one may conclude that
$\varphi$
acts trivially on
$\Gamma^\vee / 2^e \Gamma^\vee$,
as well. This implies the divisibility
$\nu_2(D_{ij}) \geq j - i + e$,
for~$i \neq j$.
\end{rem}
\subsubsection*{Third step of the proof -- A particular case}

\begin{prop}
\label{sum_of_two}
Let\/
$(\Gamma_0, b_0)$
and\/
$(\Gamma_1, b_1)$
be regular\/
\mbox{$\bbZ_2$-lattices}.
Equip the direct sum\/
$\Gamma := \Gamma_0 \oplus \Gamma_1$
with the symmetric bilinear form\/
$b := b_0 \oplus 2 b_1$.
Then~Theorem~\ref{trace_div} holds
for\/~$\Gamma$.\smallskip

\noindent
{\bf Proof.}
{\em
Let~us first note that
the~\mbox{$\smash{\bbZ_2[\sqrt{2}]}$-lattice}
\begin{align}
\label{w2_reg}
\widetilde\Gamma &:= \textstyle \left(\Gamma_0\!\otimes_{\bbZ_2}\! \bbZ_2[\sqrt{2}]\right) \oplus \frac12\sqrt{2} \left(\Gamma_1 \!\otimes_{\bbZ_2}\! \bbZ_2[\sqrt{2}]\right) \\
 &\phantom{:}\subset \Gamma \!\otimes_{\bbZ_2}\! \bbQ_2[\sqrt{2}] \,, \nonumber
\end{align}
equipped with the bilinear form induced
by~$b$,
is~regular.\smallskip

\noindent
a)
By~assumption,
$\varphi,\psi\colon \Gamma \to \Gamma$
are orthogonal maps acting trivially on
$\Gamma/2\Gamma$
and
$\Gamma/4\Gamma$,
respectively. Then~Lemma~\ref{comp_w2} shows that the induced maps
$\smash{\widetilde\varphi,\wpsi\colon \widetilde\Gamma \to \widetilde\Gamma}$
are again orthogonal, that
$(\widetilde\varphi \bmod 2)$
acts trivially on
$\smash{\widetilde\Gamma/2\widetilde\Gamma}$,
and that
$(\wpsi \bmod 2)$
acts trivially
on~$\smash{\widetilde\Gamma/4\widetilde\Gamma}$.
By~Proposition~\ref{Z2_var}.a), this implies~that
$$\smash{\Tr((\widetilde\varphi-\id) \!\circ\! (\wpsi-\id)) \in 16\bbZ_2[\sqrt{2}] \,.}$$
Since~$\smash{\widetilde\varphi}$
and~$\smash{\wpsi}$
are obtained from
$\varphi$
and~$\psi$
only by base extension, one has
$\smash{\Tr((\varphi-\id) \!\circ\! (\psi-\id)) \in 16\bbZ_2[\sqrt{2}]}$,
too. But~the latter trace is automatically
in~$\bbZ_2$,
so that the assertion~follows.\smallskip

\noindent
The~proof of~b) works along the same~lines.
}
\eop
\end{prop}

\begin{lem}
\label{comp_w2}
Let\/
$(\Gamma_0, b_0)$
and\/
$(\Gamma_1, b_1)$
be regular\/
\mbox{$\bbZ_2$-lattices}.
Equip\/~$\Gamma := \Gamma_0 \oplus \Gamma_1$
with the bilinear form\/
$b := b_0 \oplus 2 b_1$
and
let\/~$\smash{\widetilde\Gamma}$
be the\/
\mbox{$\smash{\bbZ_2[\sqrt{2}]}$-lattice,}
defined by~(\ref{w2_reg}).
Moreover,~let\/
$\varphi\colon\Gamma \to \Gamma$
be an orthogonal map, and
write\/~$\smash{\widetilde\varphi\colon \widetilde\Gamma \to \Gamma \!\otimes_{\bbZ_2}\! \bbQ_2[\sqrt{2}]}$
for the map induced
by\/~$\varphi$.

\begin{abc}
\item
Then\/~$\smash{\widetilde\varphi}$
actually sends\/
$\smash{\widetilde\Gamma}$
onto~itself.
\item
Moreover,
$\smash{\widetilde\varphi\colon \widetilde\Gamma \to \widetilde\Gamma}$
is an orthogonal~map.
\item
Let\/~$e\in\bbN$.
If\/~$(\varphi \bmod 2^e)$
acts trivially on\/
$\Gamma / 2^e \Gamma$
then\/
$\smash{(\widetilde\varphi \bmod 2^e)}$
acts trivially
on\/~$\smash{\widetilde\Gamma / 2^e \widetilde\Gamma}$.
\end{abc}\smallskip

\noindent
{\bf Proof.}
{\em
b)
is clear from the construction
of~$\smash{\widetilde\varphi}$,
once the assertion of~a) is~established.\smallskip

\noindent
a) and~c)
Let~$\calB$
be a basis
of~$\Gamma$
as in Lemma~\ref{Block_Bewertung} and~put
$\smash{(\atop{C_{00} \,C_{01}}{C_{10} \,C_{11}}) := \Mat^\calB_\calB(\varphi-\id)}$.
The~blocks
$C_{ij}$
are then matrices with coefficients
in~$\bbZ_2$.
One~even has
$\nu_2(C_{ij}) \geq e$,
for~$0 \leq i,j \leq 1$,
under the assumption of~c).

Multiplying~the basis vectors
of~$\Gamma_1$
by~$\smash{\frac12\sqrt{2}}$,
one finds a basis
$\smash{\widetilde\calB}$
of~$\smash{\widetilde\Gamma}$,
for~which
$$\Mat^{\widetilde\calB}_{\widetilde\calB}(\widetilde\varphi-\id) = \left(
\begin{array}{rr}
         C_{00} & \frac12\sqrt{2} C_{01} \\
\sqrt{2} C_{10} &                 C_{11}
\end{array}
\right) \,.$$
Thus,~in order to prove the assertions, only
$\nu_2(C_{01}) \geq 1$
and
$\nu_2(C_{01}) \geq e+1$,
respectively, need to be~verified. Both~claims are true, due to Lemma~\ref{Block_Bewertung}.
}
\eop
\end{lem}

\subsubsection*{Completion of the proof}

\begin{lem}[Change of lattice]
\label{first_reduction}
Let\/~$\Gamma$
be an arbitrary\/
\mbox{$\bbZ_2$-lattice}
and\/
$\smash{\Gamma = \bigoplus\limits_{i = 0}^N \Gamma_i}$
a decomposition as~above. Furthermore,~let\/
$\varphi\colon \Gamma \to \Gamma$
be an orthogonal~map.\smallskip

\begin{abc}
\item
Then\/
$\varphi$
maps the lattice\/
$\smash{\Gamma' := \bigoplus\limits_{i=0}^N \Gamma'_i}$,
for
$$
\Gamma'_i := 2^{-\left\lfloor \frac{i}{2}\right\rfloor} \Gamma_i \,,
$$
onto itself, as~well.
\item
Let\/~$e\in\bbN$
and suppose that\/
$\varphi$
acts trivially\/
on~$\Gamma/2^e\Gamma$.
Then\/~$\varphi$
acts trivially on the quotient\/
$\Gamma' / 2^e \Gamma'$.
\end{abc}\smallskip

\noindent\looseness-1
{\bf Proof.}
{\em
Let
$\calB$
be a basis
of~$\Gamma$
as in Lemma~\ref{Block_Bewertung}. Then~there is a basis
$\smash{\calB'}$
of
$\Gamma'$
given by scaling the basis vectors
of~$\Gamma_i$
by~$\smash{2^{-\left\lfloor \frac{i}{2}\right\rfloor}}$,
for~$i = 0, \ldots, N$.
The~matrix
$$\smash{C' = (C'_{ij})_{i,j = 0,\ldots,N} := \Mat^{\calB'}_{\calB'}(\varphi-\id)}$$
is then constructed out of the matrix
$\smash{(C_{ij})_{i,j = 0,\ldots,N} = \Mat^\calB_\calB(\varphi-\id)}$
by putting
$\smash{C'_{ij} := 2^{\left\lfloor \frac{i}{2}\right\rfloor-\left\lfloor \frac{j}{2}\right\rfloor} C_{ij}}$,
for~$i,j = 0, \ldots, N$.
This~translates the claims into certain inequalities
for~$\nu_2(C_{ij})$.\smallskip

\noindent
a)
The~assertion means that
$\smash{\nu_2(C'_{ij}) \geq 0}$,
for~$i,j = 0, \ldots, N$,
which is equivalent to
$\smash{\nu_2(C_{ij}) \geq \left\lfloor \frac{j}{2}\right\rfloor- \left\lfloor \frac{i}{2}\right\rfloor}$.
By~assumption, one has
$\smash{\nu_2(C_{ij}) \geq 0}$,
so that,
for~$i \geq j$,
there is nothing left to be~shown.
For~$i < j$,
the inequality
$\smash{\nu_2(C_{ij}) \geq \left\lfloor \frac{j}{2}\right\rfloor- \left\lfloor \frac{i}{2}\right\rfloor}$
indeed holds, due to Lemma~\ref{Block_Bewertung}.a).\smallskip

\noindent
b)
Here,~the assertion means
$\smash{\nu_2(C'_{ij}) \geq e}$,
for~$i,j = 0, \ldots, N$,
and is equivalent to
$\smash{\nu_2(C_{ij}) \geq e + \left\lfloor \frac{j}2 \right\rfloor - \left\lfloor \frac{i}2 \right\rfloor}$.
As~$\smash{\nu_2(C_{ij}) \geq e}$
holds by assumption, the case that
$i \geq j$
does not need any further~consideration. Moreover,
for~$i < j$,
the assertion follows from Lemma~\ref{Block_Bewertung}.b).
}
\eop
\end{lem}

\noindent
{\bf Proof} of {\bf Theorem~\ref{trace_div}.}
Let~$\Gamma$
be any
\mbox{$\bbZ_2$-lattice}.
In~order to prove that Theorem~\ref{trace_div} holds
for~$\Gamma$,
by~Lemma~\ref{first_reduction}, it suffices to show that exactly the same statement holds for the
\mbox{$\bbZ_2$-lattice}~$\Gamma'$.
To~actually
treat~$\Gamma'$,
note that, by construction, the bilinear
form~$\smash{b'}$
on~$\Gamma'$
is of the form
$\smash{b'_0 \oplus 2b'_1 \oplus b'_2 \oplus 2b'_2 \oplus \cdots \oplus 2^{N-2\lfloor\frac{N}2\rfloor}b'_N}$,
for regular
\mbox{$\bbZ_2$-lattices}
$\smash{(\Gamma'_0, b'_0)}$,
\ldots,
$\smash{(\Gamma'_N, b'_N)}$.
Thus,~$\Gamma'$
allows a decomposition
$\Gamma' = \Gamma''_0 \oplus \Gamma''_1$
with~$b' = b''_0 \oplus 2b''_1$,
for regular lattices
$(\Gamma''_0, b''_0)$
and~$(\Gamma''_1, b''_1)$.
But~for exactly this particular case, the assertion of~Theorem~\ref{trace_div} is established by Proposition~\ref{sum_of_two}.
\eop

\section{The point counting algorithm}
\label{Algorithm}

\begin{inp}
\begin{iii}
\item
Let a
scheme~$S'$
be given that is presented as a double cover
of~$\Pb^2_\bbZ$
branched over the union of six lines, each of which is defined
over~$\Spec\bbZ$.
Suppose~that no three of these lines have a
\mbox{$\bbQ$-rational}
point in~common, and let
$S$
be the blowing-up
of~$S'$
in the Zariski closure, equipped with the induced reduced scheme structure, of the singular locus of the generic
fibre~$\smash{S'_\bbQ}$.

Suppose that
$\smash{\Pic(S_{\overline\bbQ})}$
is a trivial
\mbox{$\smash{\Gal(\overline\bbQ/\bbQ)}$-module}
in the case that
$\smash{\rk\Pic(S_{\overline\bbQ}) \geq 16}$,
as~well.
\item
Moreover, let a
bound~$B$
be~given.
\end{iii}\smallskip

\noindent
Then~the algorithm below computes
$\#S(\bbF_{\!p})$,
for all good primes
$p < B$
of~$S$.
\end{inp}

\begin{algo}
\label{2adic_count}
\begin{ABC}
\item
{\em Initialisation.}
\begin{iii}
\item
\label{badprimes}
Calculate the odd primes
$p_1, \ldots, p_b$,
at which
$S$~has
bad~reduction. In~terms of these, declare the~map
$$\textstyle \varrho\colon \bbP \!\setminus\! \{2,p_1,\ldots,p_b\} \rightarrow \{\pm 1\}^{b+2} \,, \quad
p \mapsto ((\frac{-1}{p}), (\frac{2}{p}), (\frac{p_1}{p}), \ldots, (\frac{p_b}{p})) \,.$$
%
\item
For each
$\sigma \in \{\pm 1\}^{b+2}$,
run through
$\bbP \!\setminus\! \{2,p_1,\ldots,p_b\}$
from below, until a
prime~$l_\sigma$
is found of the kind that
$\smash{\varrho(l_\sigma) = \sigma}$.
Then~count
$\#S(\bbF_{\!l_\sigma})$
by a naive~method. From~the point~count, derive
$\smash{(\Tr(\Frob_{l_\sigma}\colon \T(S_{\overline\bbQ}, \bbZ_2) \to \T(S_{\overline\bbQ}, \bbZ_2)) \bmod 16)}$
using the Lefschetz trace formula~(\ref{Lefschetz3}) and store this value in a~table.
\end{iii}
\end{ABC}
\end{algo}

\begin{rem}
As~the action of
$\smash{\Gal(\overline\bbQ/\bbQ)}$
on
$\smash{\Pic(S_{\overline\bbQ}) \!\otimes_\bbZ\! \bbQ}$
is known to be trivial, formula~(\ref{Lefschetz3}) allows, of course, to calculate
$\smash{\Tr(\Frob_{l_\sigma}\colon \T(S_{\overline\bbQ}, \bbZ_2) \to \T(S_{\overline\bbQ}, \bbZ_2))}$
from
$\#S(\bbF_{\!l_\sigma})$~exactly.
But~only the residue class
modulo~$16$
is of~importance. Which~is exactly the information that is~stored.

Indeed,~let
$p$
be a possibly large~prime. Then,~the value
$\smash{\sigma = \varrho(p) \in \{\pm 1\}^{b+2}}$
completely determines the action of
$\Frob_p$
on
$\smash{\T(S_{\overline\bbQ}, \bbZ_2)/4\T(S_{\overline\bbQ}, \bbZ_2)}$,
according to Corollary~\ref{Fall_Q}. I.e.,~the action of
$\Frob_p$
coincides with that
of~$\Frob_{l_\sigma}$.
Furthermore, by~Theorem~\ref{Overdet}.a), this is enough to fix the trace
modulo~$16$
on~$\smash{\T(S_{\overline\bbQ}, \bbZ_2)}$.
I.e.,
$$\Tr(\Frob_p\colon \!\T(S_{\overline\bbQ}, \bbZ_2) \!\to\! \T(S_{\overline\bbQ}, \bbZ_2)) \equiv \Tr(\Frob_{l_\sigma}\colon \!\T(S_{\overline\bbQ}, \bbZ_2) \!\to\! \T(S_{\overline\bbQ}, \bbZ_2)) \pmodulo{16} \,,$$
the residue class on the right hand side being the one that was~stored. Cf.~part~\ref{pc} of the~algorithm.
\end{rem}

\begin{rem}
The~calculation of the bad primes in step~A.\ref{badprimes} involves the factorisation of a~discriminant. A~failure in this step would prevent the algorithm from~proceeding. In~our present implementation, this does not present any difficulty, as the discriminant for the family of double covers
of~$\Pb^2$
branched over six lines is highly reducible~\cite[Def.\ 7.7 and Lemma~7.8]{Yo}, cf.\ Remark~\ref{main_prop}.b). Step~A.\ref{badprimes} might, however, become an issue when trying to carry over the algorithm to other types of~surfaces. This~may concern other families of
$K3$~surfaces,
already.
\end{rem}

\begin{ttt}
The initialisation as described above ignores the group structure. It~determines
\mbox{$\smash{(\Tr(\sigma\colon \T(S_{\overline\bbQ}, \bbZ_2) \to \T(S_{\overline\bbQ}, \bbZ_2)) \bmod 16)}$}
individually, for every element of
$\Gal(K/\bbQ)$,
where~$\smash{K := \bbQ(\sqrt{-1},\sqrt{2},\sqrt{\mathstrut p_1},\ldots,\sqrt{\mathstrut p_b})}$.
The~number of these elements is exponential in the
number~$b$
of odd bad~primes.

A more efficient approach
is as follows. To~simplify notation, put
$p_{-1} := -1$
and~$p_0 := 2$.
Moreover,~let
$\sigma_{-1}, \sigma_0, \sigma_1, \ldots, \sigma_b \in \Gal(K/\bbQ)$
be the standard generators. I.e.,
$$
\textstyle \sigma_i(\sqrt{\mathstrut p_j}) =
\left\{
\begin{array}{rl}
 \!\sqrt{\mathstrut p_j} \,, & \quad\mbox{ if } j \neq i, \;j \in \{-1,0,1,\ldots,b\} \,, \\
\!-\sqrt{\mathstrut p_j} \,, & \quad\mbox{ if } j = i \,,
\end{array}
\right.
$$
for
$i \in \{-1,0,1,\ldots,b\}$.
Every~element
$\sigma \in \Gal(K/\bbQ)$
may then uniquely be described by some sequence
$s \in \{0,1\}^{b+2}$,
indexed from
$(-1)$
to~$b$,
$$\smash{\sigma = \sigma_{-1}^{s_{-1}} \sigma_0^{s_0} \sigma_1^{s_1} \cdots \sigma_b^{s_b} \,.}$$
Or,~$\sigma = \prod_{i\in M_s} \!\!\sigma_i$,
for~$M_s := \{i \in \{-1, \ldots, b\} \mid s_i = 1\}$.

With~respect to a basis of
$\smash{\T(S_{\overline\bbQ}, \bbZ_2)}$,
each generator
$\smash{\sigma_i \in \Gal(K/\bbQ)}$
yields a matrix
$E_n + 2A_i$,
with
$A_i \in \Mat_{n \times n}(\bbZ_2)$,
encoding the action on
$\smash{\T(S_{\overline\bbQ}, \bbZ_2)}$.
A~product
$\sigma_{-1}^{s_{-1}} \cdots \sigma_b^{s_b}$
then corresponds to the matrix
$$(E_n + 2 A_{-1})^{s_{-1}} \cdots (E_n + 2 A_b)^{s_b} \,,$$
the trace of which
modulo~$16$
is given by
\begin{align}
\Tr((E_n + 2 A_{-1})^{s_{-1}} \cdots (E_n + 2 A_b)^{s_b}) \equiv{}&
n + 2 \sum_{i \in M_s} \Tr(A_i) + 4\!\! \sum_{\atop{i,j \in M_s,}{i < i'}} \!\!\Tr(A_i A_{i'}) \nonumber \\[-2mm]
\label{Tr_advanced}
 & \hspace{.5cm} {}+ 8\!\!\! \sum_{\atop{i,i',i'' \in M_s,}{i < i' < i''}} \!\!\!\Tr(A_i A_{i'} A_{i''}) ~~~~~~~~~~~~~~~~~~~~~~~~~~~~ \pmodulo{16} \,.
\end{align}
Based~on~(\ref{Tr_advanced}), the traces
modulo~$8$
of all matrices
$A_i$,
together with the traces
modulo~$4$
of all products
$A_i A_{i'}$
and the traces
modulo~$2$
of all triple products
$A_i A_{i'} A_{i''}$,
can be determined efficiently by solving a system of linear congruences. This~can be described as an algorithm as~follows.
\end{ttt}

\begin{algoo}[continued]
\begin{ABC}
\addtocounter{ABC}{1}
\item[A') ]
{\em Initialisation. A more efficient approach -- First step.}

Let~$l$
run through
$\bbP \!\setminus\! \{2,p_1,\ldots,p_b\}$
from~below. Each~time, do the following.
\begin{iii}
\item
Compute~$\smash{\varrho(l)}$.
I.e.,~determine a presentation of
$\Frob_l \in \Gal(K/\bbQ)$
as a product
$\sigma_{-1}^{s_{-1}} \cdots \sigma_b^{s_b}$
of some of the standard~generators.
\item
Count~$\#S(\bbF_{\!l})$
by a naive~method. Derive
\begin{align*}
(\Tr((E_n + 2 A_{-1})^{s_{-1}} \cdots (E_n + 2 A_b))^{s_b} \bmod 16) &= \\[-1.5mm] (\Tr(\Frob_l&\colon\T(S_{\overline\bbQ}, \bbZ_2) \to \T(S_{\overline\bbQ}, \bbZ_2)) \bmod 16)
\end{align*}
from this value using the Lefschetz trace formula~(\ref{Lefschetz3}).
\item
According~to formula~(\ref{Tr_advanced}) above, write down a linear congruence involving all
$\smash{\Tr(A_i)}$,
$\smash{\Tr(A_i A_{i'})}$,
and
$\smash{\Tr(A_i A_{i'} A_{i''})}$,
for
$\smash{-1 \leq i \leq b}$,
$\smash{-1 \leq i < i' \leq b}$,
and
$\smash{-1 \leq i < i' < i'' \leq b}$,
respectively. Add~this congruence to the system of congruences already~obtained.
\item
Check~whether the system of linear congruences obtained for
$(\Tr(A_i) \bmod 8)$,
$(\Tr(A_i A_{i'}) \bmod 4)$,
and
$(\Tr(A_i A_{i'} A_{i''}) \bmod 2)$,
for
$\smash{-1 \leq i \leq b}$,
$\smash{-1 \leq i < i' \leq b}$,
and
$\smash{-1 \leq i < i' < i'' \leq b}$,
respectively, is uniquely~solvable. If~this is the case then compute the solution, store it, and terminate Step~A').
\end{iii}\smallskip

\item[A'') ]
{\em Initialisation. A more efficient approach -- Second step.}

Let~$s$
run through the elements
of~$\{0,1\}^{b+2}$.
Each~time, determine
$$\T(s) := (\Tr((E_n + 2 A_{-1})^{s_{-1}} \cdots (E_n + 2 A_b)^{s_b}) \bmod 16)$$
using formula (\ref{Tr_advanced}) and the stored values of
$(\Tr(A_i) \bmod 8)$,
$(\Tr(A_i A_{i'}) \bmod 4)$,
and
$(\Tr(A_i A_{i'} A_{i''}) \bmod 2)$.
Store~the value in a~table.

\item
\label{pc}
{\em Point counting.}

Let~$p$
run through
$\bbP \!\setminus\! \{2,p_1,\ldots,p_b\}$
from below up
to~$B$.
Each~time, do the following.
\begin{iii}
\item
Compute
$\smash{s := \varrho(p)}$,
which means to the
determine~$\Frob_p \in \Gal(K/\bbQ)$.
\item
Look up the corresponding value
$\T(s)$
in the precomputed~table. This is~just
\mbox{$\smash{(\Tr(\Frob_p\colon \T(S_{\overline\bbQ}, \bbZ_2) \to \T(S_{\overline\bbQ}, \bbZ_2)) \bmod 16)}$}.

Applying~the Lefschetz trace formula~(\ref{Lefschetz3}), calculate
$(\#S(\bbF_{\!p}) \bmod 16)$
from this value.
\item
Use a
$p$-adic
Harvey style algorithm \cite{Ha} to compute
$\smash{(\#S(\bbF_{\!p}) \bmod p)}$.
\item
Use the Chinese remainder theorem to calculate the class
\mbox{$\smash{(\#S(\bbF_{\!p}) \bmod 16p)}$}
from
$\smash{(\#S(\bbF_{\!p}) \bmod 16)}$
and~$\smash{(\#S(\bbF_{\!p}) \bmod p)}$.
\item
Determine the unique representative of this residue class
modulo~$16p$
that is compatible with Deligne's bound~(\ref{DelWeil})
for~$\#S(\bbF_{\!p})$
and output this~number.
\end{iii}
\end{ABC}
\end{algoo}

\begin{rems}[On the assumptions made
on~$S$]
\begin{iii}
\item
The~assumptions made on the intersection points of the six lines imply that
the generic
fibre~$\smash{S_\bbQ}$
of~$S$
is nonsingular, i.e.\ a 
$K3$~surface.
Moreover,~as the lines are assumed to be defined
over~$\bbZ$,
Corollary~\ref{Gal_triv_T2T}.a) shows that the action
of~$\smash{\Gal(\overline\bbQ/\bbQ)}$
on~$\smash{\T(S_{\overline\bbQ}, \bbZ_2) / 2\T(S_{\overline\bbQ}, \bbZ_2)}$
is~trivial. The~assumption
on~$\smash{\Pic(S_{\overline\bbQ})}$
is automatically fulfilled
if~$\smash{\rk\Pic(S_{\overline\bbQ}) = 16}$.
\item
The~algorithm as described immediately carries over to other types of
$K3$~surfaces,
as soon as the actions
of~$\smash{\Gal(\overline\bbQ/\bbQ)}$
on~$\smash{\Pic(S_{\overline\bbQ})}$
and
$\smash{\T(S_{\overline\bbQ}, \bbZ_2) / 2\T(S_{\overline\bbQ}, \bbZ_2)}$
are~trivial.
\item
A~further generalisation is possible to
$K3$~surfaces,
for which the action of
$\smash{\Gal(\overline\bbQ/\bbQ)}$
on~$\smash{\T(S_{\overline\bbQ}, \bbZ_2) / 2\T(S_{\overline\bbQ}, \bbZ_2)}$
is trivial and that
on~$\smash{\Pic(S_{\overline\bbQ})}$
is explicitly~known. Indeed,~the triviality of
$\smash{\Pic(S_{\overline\bbQ})}$
is only used in the references to the Lefschetz trace formula~(\ref{Lefschetz3}).
\item
Finally,~one might want to consider the case when both
\mbox{$\smash{\Gal(\overline\bbQ/\bbQ)}$-actions},
that on
$\smash{\Pic(S_{\overline\bbQ})}$
and that
on~$\smash{\T(S_{\overline\bbQ}, \bbZ_2) / 2\T(S_{\overline\bbQ}, \bbZ_2)}$,
are nontrivial, but explicitly~known. At~least when the action
on~$\smash{\T(S_{\overline\bbQ}, \bbZ_2) / 2\T(S_{\overline\bbQ}, \bbZ_2)}$
is of
exponent~$2$,
a modification of Algorithm~\ref{2adic_count} is possible, which is based on Theorem~\ref{Overdet}.b). It~may determine 
$\Tr(\Frob_p\colon \T(S_{\overline\bbQ}, \bbZ_2) \to \T(S_{\overline\bbQ}, \bbZ_2)$
only
modulo~$8$.

The~idea is as follows. Let~the number
field~$F$
be the known splitting field
of~$\smash{\T(S_{\overline\bbQ}, \bbZ_2) / 2\T(S_{\overline\bbQ}, \bbZ_2)}$.
Then~the Galois action
on~$\smash{\T(S_{\overline\bbQ}, \bbZ_2) / 4\T(S_{\overline\bbQ}, \bbZ_2)}$
factors via
$\Gal(K/\bbQ)$,
for~$K$
the maximal abelian extension of
$F$
of
exponent~$2$,
ramified only at the primes
above~$2$
and the bad primes
of~$S$.
This~is a ray class field
over~$F$
and thus, in principle, accessible to~computation. Cf.\ Section~\ref{general} for an~example.
\end{iii}
\end{rems}

\begin{rem}[Initialisation]
Among~$S_1, \ldots, S_4$,
the~surface
$S_1$
is the one having the most bad~primes. There~are actually six bad primes
$p \neq 2$.
Thus,~the direct initialisation requires to count the points for 256~primes. As~all elements
of~$\{\pm1\}^8$
have to be hit, the largest prime to be used would
be~$21\,121$.

The~more efficient approach needs to solve a system of linear congruences in 92~variables. Thus, it requires the point count only for 92~primes. As~the elements
of~$\{\pm1\}^8$
to be hit
by~$\Frob_p$
are otherwise arbitrary, we could get by working with the primes up
to~$593$.
\end{rem}

\begin{rem}
Note~that, in our examples, we always have
$2^{b+2} \leq 256$
in comparison
to~$B=10^8$.
For~other samples with many more bad primes, so that
$2^b \gg B$,
one might want to optimise by reversing steps A'') and~\ref{pc}. I.e.,~to calculate the values
$\T(\varrho(p))$
separately for each~prime.
\end{rem}

\begin{rem}[Practical performance]
\label{prac_perf}
\begin{iii}
\item
For~each of the surfaces
$S_1, \ldots, S_4$,
running~up to
$B = 10^8$,
our implementation used about 20\,GB of memory and between 8 and 12~hours of CPU time on one core of an Intel i7-7700 processor running at 3.6\,GHz.

This~running time is completely dominated by~the
modulo~$p$
point~count, for which we ran a variant of a Harvey style algorithm in
\mbox{$p$-adic}
precision~$1$
with remainder tree implemented for this particular~project. Note~that our implementation is in {\tt magma}, not in a compiled language. And~that, presumably, some of the possible optimisations are still~missing. Cf.\ \cite[Sections 3 and~4]{EJ16} for a description of an earlier implementation. 
\item
The~initialisation, as described in A') and~A''), took less than one minute per surface. More~precisely, for the
surface~$S_1$,
the naive point counting had to be done for
$92$~primes,
which took around
$49$~seconds.
For~all other steps of the initialisation together, including the linear algebra calculations, the {\tt magma} profiler reports a running time
of
$2.5$~seconds.
For~the surfaces
$S_2$,
$S_3$,
and~$S_4$,
the initialisation runs faster by a factor of at
least~$10$,
because there are fewer bad~primes.

The actual
modulo~$16$
point counting mainly required Legendre symbol computations for the slightly more than~5.7 million primes up
to~$B = 10^8$,
which took only
$23$~seconds
per~surface. Finally,~the determination of the point counts
modulo~$16p$
took
$5$~seconds,
which are essentially accounted for the computations related to the Chinese remainder~theorem.
\end{iii}
\end{rem}

\begin{rem}[Point counting
modulo~$p^2$
versus
modulo~$p$]
\label{factor10}
\begin{iii}
\item
\mbox{$p$-adic}
point counting for a surface of the shape
$w^2 = xyz f_3(x,y,z)$
requires to do the following: In~order to count
modulo~$p$,
one has to compute the coefficient at
$(xyz)^{(p-1)/2}$
in~$\smash{f_3^{(p-1)/2}}$
with
\mbox{$p$-adic}
precision~$1$.
In all our examples, the moving simplex approach (cf.\ \cite[Remark~4.8]{EJ16}) never resulted in a
loss of
\mbox{$p$-adic}
precison. Thus, we were able to work with
\mbox{$p$-adic}
precision~$1$
during all the intermediate~steps.

On~the other hand, for point counting
modulo~$p^2$,
one has to compute the coefficient at
$(xyz)^{(p-1)/2}$
in~$\smash{f_3^{(p-1)/2}}$
with
\mbox{$p$-adic}
precision~$2$
and, furthermore, the coefficient at
$x^{i(p-1)/2} y^{j(p-1)/2} z^{k(p-1)/2}$
in~$\smash{f_3^{3(p-1)/2}}$,
for every triple
$(i,j,k) \in \bbN^3$
of odd numbers such that
$i+j+k=9$.
Thus, instead of computing one coefficient, one has to compute eleven. Assuming that this can be done without
\mbox{$p$-adic}
precision loss, one can work with
\mbox{$p$-adic}
precision~$2$
during all the intermediate~steps. This~indicates that one has to expect an increase of the run time by at least a factor
of~$22$.

As the exponent
in~$\smash{f_3^{3(p-1)/2}}$
is increased by a factor
of~$3$
compared
to~$\smash{f_3^{(p-1)/2}}$,
a naive implementation would slow down the process even more, in the worst case by another factor
of~$3$.
However, a better implementation using multipoint evaluation techniques \cite[Section~10.1]{vGG} might reduce this factor significantly.
\item
For checking correctness, we implemented a simple
\mbox{$p$-adic}
point counting with
\mbox{$p$-adic}
precision~$2$
that does not use advanced techniques such as the remainder tree. This~implementation took about one day of CPU time per surface, running only
to~$B = 10^5$.
It~is important to note, however, that this is not a fair comparison, since too many optimisations were~missing.
\end{iii}
\end{rem}

\begin{rems}[Results]
\begin{abc}
\item
The main outcome of our computations are the histograms presented in~\cite[Section~5]{EJ21}.
\item
The~distribution of the traces
modulo~$16$
relative to the elements
in~$\{\pm1\}^{b+2}$,
as indicated in the table below, appears to be rather~erratic.

\begin{table}[H]
\begin{center}
\setlength{\tabcolsep}{1.5mm}
\begin{tabular}{|c||r|r|r|r|r|r|r|r|r|r|r|r|r|r|r|r|}
\hline
Residue mod 16 & 0 & 1 & 2 & \phantom{1}3 & 4 & 5 & 6 & \phantom{1}7 & 8 & \phantom{1}9 & 10 & 11 & 12 & 13 & 14 & 15 \\\hline\hline
$S_1$ & 46 &  0 & 44 & 0 & 32 &  0 & 26 & 0 & 18 & 0 & 36 & 0 & 32 &  0 & 22 & 0 \\\hline
$S_2$ & 0 &  0 &  7 & 0 &  0 &  0 &  7 & 0 &  0 & 0 &  9 & 0 &  0 &  0 &  9 & 0 \\\hline
$S_3$ & 0 & 20 &  0 & 0 &  0 & 28 &  0 & 0 &  0 & 4 &  0 & 0 &  0 & 12 &  0 & 0 \\\hline
$S_4$ & 8 &  0 &  2 & 0 &  0 &  0 &  4 & 0 &  0 & 0 &  2 & 0 &  0 &  0 &  0 & 0 \\\hline
\end{tabular}\vskip-3mm
\caption{Number of elements of
$\{\pm1\}^{b+2}$
for each residue
$\bmod 16$}
\end{center}
\end{table}\vskip-5mm

Nevertheless,~there~are a few more observations that should perhaps be~noticed.
\begin{iii}
\item
For~instance,
for~$S_3$,
it happens that
\mbox{$(\Tr(\Frob_p\colon \!\T(S_{\overline\bbQ}, \bbZ_2) \!\to\! \T(S_{\overline\bbQ}, \bbZ_2)) \bmod 16)$}
is independent of the Legendre symbol
$\smash{(\frac{2}{p})}$.
For~$S_4$,
it suffices to consider
$\smash{(\frac{6}{p})}$,
instead of
$\smash{(\frac{2}{p})}$
and~$\smash{(\frac{3}{p})}$
individually. These,~however, are the only regularities that occurred above those predicted by Corollary~\ref{Fall_Q}.
\item
(Explanation of the zeroes in Table~1)
For~$p$
odd, a double cover of
$\smash{\Pb^2_{\bbF_{\!p}}}$,
branched over six
\mbox{$\bbF_{\!p}$-rational}
lines in general position, has an odd number of points, since the branch locus~has. This yields that
$\smash{\Tr(\Frob_p\colon \T(S_{\overline\bbQ}, \bbZ_2) \to \T(S_{\overline\bbQ}, \bbZ_2))}$
is always even in the case of geometric Picard
rank~$16$,
and odd, for geometric Picard
rank~$17$.

But~more is~true. One~has that
$\det(E_n + 2A) = \pm1$
implies
$1 + 2\Tr(A) \equiv \pm1 \pmodulo 4$
and therefore
$\Tr(E_n + 2A) \equiv n - 1 \pm1 \pmodulo 4$.
This~explains why, for
$S_1$
and~$S_4$,
there are equally many elements
of~$\{\pm1\}^{b+2}$
leading to traces
$(0 \bmod 4)$
and~\mbox{$(2 \bmod 4)$}.

The~surfaces
$S_2$~and~$S_3$,
however, have trivial jump characters~\cite{CEJ}, so that
$\smash{\det(\Frob_p\colon \T(S_{\overline\bbQ}, \bbZ_2) \to \T(S_{\overline\bbQ}, \bbZ_2)) = +1}$,
for every
prime~$p$.
This~explains why only traces
$(2 \bmod 4)$
occur
for~$S_2$
and only traces
$(1 \bmod 4)$
for~$S_3$.

Finally,~for every
$p \equiv 3 \pmodulo 4$,
$\smash{\Tr(\Frob_p\colon \T(S_{4,\overline\bbQ}, \bbZ_2) \to \T(S_{4,\overline\bbQ}, \bbZ_2)) = 0}$
holds exactly.
This~is the spike in \cite[Figure~3]{EJ21}. On~the other hand,
for~$p \equiv 1 \pmodulo 4$,
one has
$\smash{\det(\Frob_p\colon \T(S_{4,\overline\bbQ}, \bbZ_2) \to \T(S_{4,\overline\bbQ}, \bbZ_2)) = +1}$,
so that only traces
$(2 \bmod 4)$
are~allowed. The~non-occurrence of traces
$(14 \bmod 16)$
seems to be explained only by the fact that
$b=2$
is very~small.
\end{iii}
\end{abc}
\end{rems}

\section{A more advanced example -- A surface with real multiplication}
\label{general}

This section is devoted to the surface
$S_5$
from Example~\ref{histo_exs}.v).

\begin{ttt}[Properties of the surface]
\begin{iii}
\item
This~surface is \cite[Example~5.8]{EJ21}. Also,
$S_5$
is isomorphic to the specialisation to
$t=0$
of the family described in \cite[Example 1.5]{EJ20a}. In~particular,
$\smash{\rk\Pic(S_{5,\overline\bbQ}) = 16}$
and
$S_5(\bbC)$
has real multiplication
by~$\smash{\bbQ(\sqrt{5})}$.
\item
The~surface
$S'_{5,\bbQ}$
is a double cover of
$\smash{\Pb^2_\bbQ}$,
branched geometrically over six lines, any three of which do not have a geometric point in~common. Two~of these lines are defined
over~$\bbQ$,
while the other four are defined
over~$\bbQ(\zeta_5)$
and permuted cyclically
by~$\Gal(\bbQ(\zeta_5)/\bbQ) \cong \bbZ/4\bbZ$.

By~Theorem~\ref{S6_action}.a.i) and~b),
$\smash{\Pic(S_{5,\overline\bbQ})}$
contains a sublattice of full rank that is a linear permutation representation
of~$\Gal(\bbQ(\zeta_5)/\bbQ)$.
The~underlying permutation representation has one fixed point and is otherwise the action on pairs of six objects, two of which are fixed while the others form a
\mbox{$4$-cycle}.
Thus,~in total, there~are three orbits of size four each, one orbit of size two, and two fixed points. I.e.,
\begin{equation}
\label{Tr_alg}
\Tr(\sigma\colon \Pic(S_{5,\overline\bbQ}) \!\otimes_\bbZ\! \bbQ \to \Pic(S_{5,\overline\bbQ}) \!\otimes_\bbZ\! \bbQ) = \left\{\!
\begin{array}{rc}
16 \,, & ~{\rm if } \ord(\sigma) = 1 \,, \\
 4 \,, & ~{\rm if } \ord(\sigma) = 2 \,, \\
 2 \,, & ~{\rm if } \ord(\sigma) = 4 \,.
\end{array}
\right.
\end{equation}
\item
For every good prime
$p \equiv 2,3 \pmod 5$,
one~has
$$\#S_5(\bbF_{\!p}) = p^2 + 2p + 1$$
by \cite[Lemma~6.7]{EJ20a}. This~means that
$\smash{\Tr(\Frob_p\colon \T(S_{5,\overline\bbQ}, \bbZ_2) \to \T(S_{5,\overline\bbQ}, \bbZ_2)) = 0}$.
\item
In~other words, only the good primes
$p \equiv 1,4 \pmod 5$
need consideration in this~example. For~these, according to~(\ref{Tr_alg}) and~(\ref{Lefschetz3}), one~has
\begin{equation}
\label{Lefsch}
\#S_5(\bbF_{\!p}) = p^2 + T_\alg p + \Tr(\Frob_p\colon \T(S_{5,\overline\bbQ}, \bbZ_2) \to \T(S_{5,\overline\bbQ}, \bbZ_2)) p + 1 \,,
\end{equation}
for
$$
T_\alg := \left\{\!
\begin{array}{rc}
16 \,, & ~\;{\rm if }\; p \equiv 1 \pmodulo 5 \,, \\
 4 \,, & ~\;{\rm if }\; p \equiv 4 \pmodulo 5 \,.
\end{array}
\right.
$$
\end{iii}
\end{ttt}

\begin{ttt}[The
$\Gal(\overline\bbQ/\bbQ)$-action
on
$\smash{\T(S_{5,\overline\bbQ}, \bbZ_2)/2\T(S_{5,\overline\bbQ}, \bbZ_2)}$]
\leavevmode

\noindent
Theorem~\ref{S6_action}.a.i) and~b), together with Theorem~\ref{prop_Br2}, provides an explicit description of
$\smash{\T(S_{\overline\bbQ}, \bbZ_2)/2\T(S_{\overline\bbQ}, \bbZ_2)}$
as a
\mbox{$\smash{\Gal(\overline\bbQ/\bbQ)}$-module}.
The~results are as~follows.

\begin{iii}
\item
The
\mbox{$\smash{\Gal(\overline\bbQ/\bbQ)}$-action}
on~$\smash{\Br(S_{\overline\bbQ})_2}$
factors via
$\Gal(\bbQ(\zeta_5)/\bbQ)$,
which is cyclic of order~four.
\item
With~respect to a suitable basis, the action of a generator
of~$\Gal(\bbQ(\zeta_5)/\bbQ)$
on~$\smash{\Br(S_{\overline\bbQ})_2}$
is given by the~matrix
$$
\left(
\begin{smallmatrix}
1 & 1 & 1 & 0 & 1 & 1 \\
0 & 1 & 1 & 1 & 1 & 1 \\
1 & 0 & 1 & 1 & 1 & 1 \\
1 & 1 & 0 & 1 & 1 & 1 \\
1 & 1 & 1 & 1 & 0 & 1 \\
1 & 1 & 1 & 1 & 1 & 0
\end{smallmatrix}
\right) \,.
$$
Thus,~the action of the element of order two is given~by
\begin{equation}
\label{6x6}
\left(
\begin{smallmatrix}
0 & 1 & 0 & 0 & 0 & 0 \\
1 & 0 & 0 & 0 & 0 & 0 \\
0 & 0 & 0 & 1 & 0 & 0 \\
0 & 0 & 1 & 0 & 0 & 0 \\
0 & 0 & 0 & 0 & 1 & 0 \\
0 & 0 & 0 & 0 & 0 & 1
\end{smallmatrix}
\right) \,.
\end{equation}
\item
By~Theorem~\ref{prop_Br2}, one has
$\smash{\T(S_{\overline\bbQ}, \bbZ_2)/2\T(S_{\overline\bbQ}, \bbZ_2) \cong (\Br(S_{\overline\bbQ})_2)^\vee}$.
The~action~on
$\smash{\T(S_{\overline\bbQ}, \bbZ_2)/2\T(S_{\overline\bbQ}, \bbZ_2)}$
is therefore provided by the transposed inverses of the matrices given.
\end{iii}
\end{ttt}

\begin{ttt}[The
$\Gal(\overline\bbQ/\bbQ)$-action
on
$\smash{\T(S_{\overline\bbQ}, \bbZ_2)/4\T(S_{\overline\bbQ}, \bbZ_2)}$]
\leavevmode

\noindent
As the Galois action on
$\smash{\T(S_{\overline\bbQ}, \bbZ_2)/2\T(S_{\overline\bbQ}, \bbZ_2)}$
factors via
$\Gal(\bbQ(\zeta_5)/\bbQ)$,
the action on
$\smash{\T(S_{\overline\bbQ}, \bbZ_2)/4\T(S_{\overline\bbQ}, \bbZ_2)}$
factors via
$\Gal(K/\bbQ)$,
for
$K$
the maximal abelian extension of
$\bbQ(\zeta_5)$
of
exponent~$2$,
ramified only
at~$2$
and the prime
above~$5$.
It~is not hard to see~that
\begin{equation}
\label{split_field}
K = \bbQ(\sqrt{-1}, \sqrt{2}, \zeta_5, \sqrt{\zeta_5-1}, \sqrt{\smash{\zeta_5^2}-1}) \,.
\end{equation}
Note~here that the prime
above~$5$
is~$(\zeta_5-1)$.
Moreover,~the unit group
$\bbZ[\zeta_5]^*$
is generated by the cyclotomic unit
$\smash{\frac{\zeta_5^2-1}{\zeta_5-1} = \zeta_5+1 }$,
together
with~$(-\zeta_5)$.
\end{ttt}

\begin{ttt}[Adaptation of the point counting algorithm]
\leavevmode

\noindent
Suppose~a good prime
$p \equiv 1,4 \pmod 5$
to be~given. Then~one determines the conjugacy class
$\Frob_p$
in~$\Gal(K/\bbQ)$.
As~before, there is a small prime
$l$
such that
$\Frob_l$
is the same~class. By~the construction
of~$K$,
this means that the action of
$\Frob_p$
on
$\smash{\T(S_{5,\overline\bbQ}, \bbZ_2)/4\T(S_{5,\overline\bbQ}, \bbZ_2)}$
agrees with that
of~$\Frob_l$.
Moreover,~one may look up the value
$\smash{(\Tr(\Frob_l\colon \T(S_{5,\overline\bbQ}, \bbZ_2) \to \T(S_{5,\overline\bbQ}, \bbZ_2)) \bmod 16)}$
in a precomputed~table. Let~us now distinguish between the two~cases.\smallskip

\noindent
{\bf Case~1:}
$p \equiv 4 \pmod 5$.
Then
$\smash{\Frob_p \in \Gal(\bbQ(\zeta_5)/\bbQ)}$
is the element of order~two. Correspondingly,~the action of
$\smash{\Frob_p}$
on~$\smash{\T(S_{5,\overline\bbQ}, \bbZ_2)/2\T(S_{5,\overline\bbQ}, \bbZ_2)}$
is not trivial, cf.\ (\ref{6x6}), but that of
$\smash{\Frob_p^2}$~is.
Hence,~Theorem~\ref{Overdet}.b) applies and shows that
\begin{align*}
\Tr(\Frob_p\colon \T(S_{5,\overline\bbQ}, \bbZ_2) \to \T(S_{5,\overline\bbQ}, \bbZ_2)) &\equiv \\[-1.2mm]
\Tr(\Frob_l&\colon \!\T(S_{5,\overline\bbQ}, \bbZ_2) \to \T(S_{5,\overline\bbQ}, \bbZ_2)) \pmodulo 8 \,.
\end{align*}
Thus, the precomputation determines
$\smash{(\Tr(\Frob_p\!\colon \!\!\T(S_{5,\overline\bbQ}, \bbZ_2) \!\to\! \T(S_{5,\overline\bbQ}, \bbZ_2)) \!\bmod 8)}$,
and therefore
$(\#S_5(\bbF_{\!p}) \bmod 8)$,
when taking (\ref{Lefsch}) into~consideration. Combining~this with a point count
modulo~$p$,
one may easily compute
$(\#S_5(\bbF_{\!p}) \bmod 8p)$.

Recall~at this point that
$S_5(\bbC)$
has real multiplication by a quadratic number~field. This~causes the algebraic monodromy group to be significantly smaller than~usual. Concretely, the Zariski closure of the image
of~$\smash{\Gal(\overline\bbQ/\bbQ)}$
in~$\smash{\Aut(\T(S_{\overline\bbQ}, \overline\bbQ_2))}$
is isomorphic to
$\smash{[\O_3(\overline\bbQ_2)]^2}$,
and not
to~$\smash{\O_6(\overline\bbQ_2)}$.
Moreover,~as
$p \equiv 4 \pmodulo 5$,
the action
of~$\Frob_p$
lies
in~$\smash{[\O_3^-(\overline\bbQ_2)]^2}$
\cite[Theorem~5.9]{EJ21}. In~particular, two of the six eigenvalues are bound to
be~$(-1)$.

Denoting~the other eigenvalues by
$\lambda_1, \ldots, \lambda_4$,
formula~(\ref{Lefsch}) takes the~form
$$\#S_5(\bbF_{\!p}) = p^2 + (2 + \lambda_1 + \lambda_2 + \lambda_3 + \lambda_4)p + 1 \,.$$
As~$|\lambda_i| = 1$,
for~$i = 1, \ldots, 4$,
this shows that
$(\#S_5(\bbF_{\!p}) \bmod 8p)$
uniquely determines
$\#S_5(\bbF_{\!p})$,
except for the two edge cases
$\lambda_1 = \cdots = \lambda_4 = 1$
and
$\lambda_1 = \cdots = \lambda_4 = -1$,
which seem indistinguishable. However,~the second of these does not occur, due to the Lemma~below.
\end{ttt}

\begin{lem}
Let\/~$p$
be an odd prime number and\/
$S$
a\/
$K3$~surface
over\/~$\bbF_{\!p}$
as in~\ref{six_lines}. Suppose~that\/
$\Frob_p$
acts on the six lines as a permutation of order at most~two. Then,~among the six eigenvalues
of\/~$\Frob_p$
on the orthogonal complement of\/
$\pi^*[l]$,
$\pi^*[e_{12}]$,
\ldots,
$\pi^*[e_{56}]$
in\/~$\smash{H^2_\et(S_{\overline\bbF_{\!p}}, \bbZ_2(1))}$,
at least one is not equal
to\/~$(-1)$.\smallskip

\noindent
{\bf Proof.}
{\em
Assume~the contrary. Then,~for the arithmetic Picard group
of~$S$,
one has
$\Pic(S) \!\otimes_\bbZ\! \bbQ \subseteq \bbQ \pi^*[l] \oplus \bbQ \pi^*[e_{12}] \oplus\cdots\oplus \bbQ \pi^*[e_{56}]$.
An~orthogonal basis is provided by the class
$\pi^*[l]$,
the classes
$\smash{\pi^*[e_{ij}]}$
being invariant
under~$\Frob_p$,
and the classes
$\smash{\pi^*[e_{ij}] + \pi^*[e_{i'j'}]}$
formed by an orbit of size~two. The~self-intersection numbers of these are equal~to
$2$,
$(-2)$,
and~$(-4)$,
respectively, so that the discriminant
of~$\Pic(S)$
as a quadratic space
is~$\smash{(\pm\overline{1})}$
or~$\smash{(\pm\overline{2}) \in \bbQ^*/\bbQ^* {}^2}$.

On~the other hand, according to the Artin--Tate formula~\cite[Theorem~6.1]{Mi}, the value of this discriminant may be calculated as
$\smash{\overline{\pm p \!\cdot\! \prod_{\lambda\neq1} (1-\lambda)}}$,
the product being taken over all eigenvalues
$\lambda \neq 1$
of~$\Frob_p$
on~$\smash{H^2_\et(S_{\overline\bbF_{\!p}}, \bbZ_2(1))}$,
counted with multiplicities. However,~according to our assumptions, every eigenvalue
$\neq \!1$
of~$\Frob_p$
is~$(-1)$,
which enforces the discriminant
of~$\Pic(S)$
to be
$\smash{(\pm\overline{p})}$
or~$\smash{(\pm\overline{2p})}$.
A~contradiction.%
}%
\eop
\end{lem}

\noindent
{\bf Case~2:}
$p \equiv 1 \pmod 5$.
This case is easier.
One has that
$\Frob_p \in \Gal(\bbQ(\zeta_5)/\bbQ)$
is the neutral~element. Consequently,~the action of
$\Frob_p$
on~$\smash{\T(S_{5,\overline\bbQ}, \bbZ_2)/2\T(S_{5,\overline\bbQ}, \bbZ_2)}$
is trivial, so that Theorem~\ref{Overdet}.a) applies and~shows
\begin{align*}
\Tr(\Frob_p\colon \T(S_{5,\overline\bbQ}, \bbZ_2) \to \T(S_{5,\overline\bbQ}, \bbZ_2)) \equiv& \\[-1.2mm]
\Tr(\Frob_l&\colon \T(S_{5,\overline\bbQ}, \bbZ_2) \to \T(5,S_{\overline\bbQ}, \bbZ_2)) \pmodulo {16} \,.
\end{align*}
I.e.,~the precomputed value fixes
$\smash{(\Tr(\Frob_p\colon \T(S_{5,\overline\bbQ}, \bbZ_2) \to \T(S_{5,\overline\bbQ}, \bbZ_2)) \bmod 16)}$
and therefore
$(\#S_5(\bbF_{\!p}) \bmod 16)$.
Combining~this with a
modulo~$p$
point count, one may compute
$(\#S_5(\bbF_{\!p}) \bmod 16p)$.
And,~similarly to (\ref{DelWeil}), this is enough to completely determine
$\#S_5(\bbF_{\!p})$.
Note~that, in formula~(\ref{Lefsch}), the trace
of~$\Frob_p$
is bounded
by~$6$
in absolute~value.

\begin{rem}
The~information on
$(\#S_5(\bbF_{\!p}) \bmod 8p)$
suffices to determine
$\#S_5(\bbF_{\!p})$
in Case~2, as~well. Indeed, real multiplication causes two further
eigenvalues~$(+1)$,
so that one has
$\#S_5(\bbF_{\!p}) = p^2 + (18 + \lambda_1 + \lambda_2 + \lambda_3 + \lambda_4)p + 1$.
This~shows that, again, only the two edge cases seem indistinguishable.
But~$22$
eigenvalues~$(+1)$
are impossible for a
$K3$~surface
over the prime
field~$\bbF_{\!p}$,
for
$p \neq 2$,
\cite[(6.8)]{Ar}.
\end{rem}

\begin{rem}[Practical performance]
The~total running time for
$S_5$
was around 58~hours, which is a lot more than for the other examples. Cf.\ Remark~\ref{prac_perf}.

The~difference comes mainly from the
\mbox{$p$-adic}
algorithm. In~fact, for a surface given by
$W^2 = f_6(T_1,T_2,T_3)$
as a double cover
of~$\smash{\Pb^2_{\bbF_{\!p}}}$,
the number of
\mbox{$\bbF_{\!p}$-rational}
points
modulo~$p$
depends only on the coefficient at
$\smash{T_1^{p-1}T_2^{p-1}T_3^{p-1}}$
of~$\smash{f_6^{(p-1)/2}}$.
Therefore,\vspace{.3mm}
for~$f_6 = T_1T_2T_3 f_3$,
one only needs to compute the coefficient at
$\smash{T_1^{(p-1)/2}T_2^{(p-1)/2}T_3^{(p-1)/2}}$
of~$\smash{f_3^{(p-1)/2}}$.
Our implementation makes systematic use of this simplification, which applies to
$S_1$,
$S_2$,
$S_3$,
and~$S_4$.

However,~the equation\vspace{.3mm}
of~$S_5$
only has the form
$f_6 = T_1T_2 f_4$,
so that the coefficient at
$\smash{T_1^{(p-1)/2} T_2^{(p-1)/2} T_3^{p-1}}$
of
$\smash{f_4^{(p-1)/2}}$
is asked~for. This~computation is more elaborate, so the last example took about five times~longer.
\end{rem}

\begin{rems}[Results]
\begin{abc}
\item
The~main purpose of our computations for this example was to generate the histogram in~\cite[Figure~4]{EJ21} to the~left.
\item
Moreover,~our computations show that, in this particular example, the trace
$\smash{\Tr(\Frob_p\colon \T(S_{5,\overline\bbQ}, \bbZ_2) \to \T(S_{5,\overline\bbQ}, \bbZ_2))}$
modulo~$8$
or~$16$,
respectively, is determined by the conjugacy class
$\Frob_p$
in a field a lot smaller than the
field~$K$
deduced from the general theory, cf.\ formula~(\ref{split_field}). In~fact, the following~holds.
\begin{iii}
\item
If~$p \equiv 4 \pmodulo 5$
then
\begin{align*}
\Tr(\Frob_p\colon \T(S_{5,\overline\bbQ}, \bbZ_2) \to \T(S_{5,\overline\bbQ}, \bbZ_2)) &\equiv \\
&\left\{
\begin{array}{rr}
6 \pmodulo 8 \,, & {\rm \;\;if\;} (-1) {\rm \;is\;a\;square\;in\;} \bbF_{\!p} \,,\\
2 \pmodulo 8 \,, & {\rm otherwise} \,.
\end{array}
\right.
\end{align*}
\item
If~$p \equiv 1 \pmodulo 5$
then
\begin{align*}
\;\;\Tr&(\Frob_p\colon \T(S_{5,\overline\bbQ}, \bbZ_2) \to \T(S_{5,\overline\bbQ}, \bbZ_2)) \equiv \\
&\left\{
\begin{array}{rr}
6 \pmodulo {16} \,, & {\rm \;if\;} (-1), \;(\zeta_5-1), {\rm \;and\;} (\zeta_5^2-1) {\rm \;are\;squares\;in\;} \bbF_{\!p} \,,\\
2 \pmodulo {16} \,, & {\rm \;\;if\;} (-1) {\rm \;is\;a\;square,\,but\;} {\smash{\frac{\zeta_5^2-1}{\zeta_5-1}} = \zeta_5+1} {\rm\text{\;is\;a\;non-square\;in\;}} \bbF_{\!p} \,,\\
14 \pmodulo {16} \,, & {\rm otherwise} \,.
\end{array}
\right.
\end{align*}
\end{iii}
Note~that these conditions are independent of the choice of the fifth root of unity
$\zeta_5 \in \bbF_{\!p}$.
Indeed,~replacing
$\zeta_5$
by~$\zeta_5^2$,
one finds that only
$\smash{\frac{\zeta_5^4-1}{\zeta_5-1} = -\zeta^4 = -(\zeta^2)^2}$
needs to be identified as being a~square.
\end{abc}
\end{rems}

\section{The very general picture}
\label{very_general}

Consider the class of all
$K3$~surfaces
over~$\bbQ$,
or even a different kind of surfaces, but suppose that
$H^1_\et(S_{\overline\bbQ}, \bbZ_2) = 0$
and that
$\Pic(S_{\overline\bbQ})$
is computable as a
$\Gal(\overline\bbQ/\bbQ)$-module.
Is~it then, at least in principle, possible to devise a
\mbox{$2$-adic}
point counting algorithm for such a class of~surfaces? This,~in essence, means to make the
$\Gal(\overline\bbQ/\bbQ)$-module
$\smash{\T(S_{\overline\bbQ}, \bbZ_2) / 2^m \T(S_{\overline\bbQ}, \bbZ_2)}$
explicit, for a suitable value
of~$m$.

A~major portion of the information on the 
$\Gal(\overline\bbQ/\bbQ)$-module
structure is encoded in the splitting
field~$F_m$
of~$\T(S_{\overline\bbQ}, \bbZ_2) / 2^m \T(S_{\overline\bbQ}, \bbZ_2)$.
This~is, by definition, the smallest field allowing a commutative~diagram
$$
\diagram
\Gal(\overline\bbQ/\bbQ) \rrto^{\varrho\;\;\;\;\;\;\;\;\;\;\;\;\;\;\;\;\;\;\;\;\;} \ar@{->>}[d] && \Aut(\T(S_{\overline\bbQ}, \bbZ_2) / 2^m \T(S_{\overline\bbQ}, \bbZ_2)) \\
\Gal(F_m/\bbQ) \urrto^{\varrho_{F_m}\;\;\;\;\;\;}
\enddiagram
\,,
$$
for
$\varrho$
the natural~action. One~has that
$F_m$
is an algebraic number field, for every
$n\in\bbN$.

\begin{ttt}
One~might want to determine the splitting
fields~$F_m$
inductively. The~induction step from
$m$
to~$m+1$
should work as~follows.

One~has that
$\smash{\{A \in \GL_n(\bbZ / 2^{m+1} \bbZ) \mid\! A \equiv E_n\!\! \pmodulo {2^m}\}}$
is an elementary abelian
\mbox{$2$-group},
hence
$F_{m+1}/F_m$
is always an abelian field extension of
exponent~$2$.
Moreover,~the smooth specialisation theorem~\cite[Expos\'e XVI, Corollaire~2.3]{SGA4} implies that
$F_{m+1}/F_m$
is unramified at any prime of odd residue characteristic, except possibly those lying above the prime numbers at which
$S$~has
bad~reduction. In~other words, an upper bound
for~$F_{m+1}$
is provided by a certain ray class field
of~$F_m$,
which is, at least in principle, amenable to~computation.
\end{ttt}

\begin{ttt}
On~the other hand, the splitting
field~$F_2$
of~$\smash{\T(S_{\overline\bbQ}, \bbZ_2) / 2 \T(S_{\overline\bbQ}, \bbZ_2)}$
is a number field of degree
$\leq \!\#\GL_t(\bbF_{\!2}) = 2^{t(t-1)/2} (2^t-1)\cdots(2^1-1)$,
for~$t := \dim T$,
unramified at every odd prime of good reduction
of~$S$.
There~are only finitely many such number fields, according to Minkowski's Theorem~\cite[Theorem III.2.13]{Ne}, and to determine all of them is, in theory, effective. Thus,~the composite of all these fields is an upper bound
for~$F_2$.

Such~an approach, however, appears practically unfeasable under virtually all circumstances. Thus,~it seems that, generally speaking, the base case is more complex than the induction~step.
\end{ttt}

\begin{rem}
In~order to settle this issue with the base case for the particular family of
$K3$~surfaces
considered in this article, we decided to apply the isomorphism from Theorem~\ref{prop_Br2}. This~requires to make
$\Br(S_{\overline\bbQ})_2$
explicit, for which there is no obvious general approach~either. The~work of A.\,N.~Skorobogatov \cite[Theorem~1.1]{Sk} we use is limited to double~covers. Furthermore,~it provides, in general, only an exact sequence, which might be non-split in certain~cases.
\end{rem}

\begin{rem}
Only a subfield
of~$F_m$,
the {\em trace field,} is relevant for the~algorithm. This~is a minimal
field~$K$,
for which there is a commutative~diagram
$$
\diagram
\Gal(F_m/\bbQ) \rrto^{\varrho_{F_m}\;\;\;\;\;\;\;\;\;\;\;\;\;\;\;\;\;\;} \ar@{->>}[d] && \Aut(\T(S_{\overline\bbQ}, \bbZ_2) / 2^m \T(S_{\overline\bbQ}, \bbZ_2)) \dto^\Tr \\
\Gal(K/\bbQ) \rrto && \bbZ/2^m\bbZ
\enddiagram
\,.
$$
The
\mbox{$2$-adic}
overdetermination phenomenon established in Section~\ref{overdet} indicates that
$[K:\bbQ]$
may be significantly smaller
than~$[F_m:\bbQ]$.
\end{rem}

\setlength\parindent{0mm}

\frenchspacing
\begin{thebibliography}{SGA4\textonehalf}
\bibitem[Ar]{Ar}
Artin, M.: {\em Supersingular
$K3$
surfaces,} Ann.\ Sci.\ \'Ecole Norm.\ Sup.\ {\bf 7}\br(1974)\brr543--567
\bibitem[BCP]{BCP}
Bosma, W., Cannon, J., and Playoust, C.: {\em The Magma algebra system~I.\ The user language,} J.\ Symbolic Comput.\ {\bf 24}\br(1997)\brr235--265
\bibitem[Ch]{Ch}
Charles, F.: {\em The Tate conjecture for
$K3$~surfaces
over finite fields,} Invent.\ Math.\ {\bf 194}\br(2013)\brr119--145
\bibitem[CS]{CS}
Conway, J.\,H. and Sloane, N.\,J.\,A.: Sphere packings, lattices and groups, Third edition, Grundlehren der Mathematischen Wissenschaften 290, {\em Springer,} New York~1999
\bibitem[CEJ]{CEJ}
Costa, E., Elsenhans, A.-S., and Jahnel, J.: {\em On the distribution of the Picard ranks of the reductions of a
$K3$~surface},
Research in Number Theory {\bf 6}\br(2020)art.27
\bibitem[De74]{De74}
Deligne,~P.: {\em La conjecture de Weil~I,} Publ.\ Math.\ IHES {\bf 43}\br(1974)\brr273--307
\bibitem[De81]{De81}
Deligne, P.: {\em Rel\`evement des surfaces\/
$K3$
en caract\'eristique nulle,} Prepared for publication by Luc Illusie, in: Algebraic surfaces (Orsay 1976--78), Lecture Notes in Math.~868, {\em Springer,} Berlin--New York~1981, 58--79
\bibitem[DM]{DM}
Dixon, J.\,D.\ and Mortimer, B.: Permutation groups, Graduate Texts in Mathematics~163, {\em Springer,} New York~1996
\bibitem[EJ16]{EJ16}
Elsenhans, A.-S.\ and Jahnel, J.: {\em Point counting on
$K3$~surfaces
and an application concerning real and complex multiplication,} in: Proceedings of the ANTS XII conference (Kaiserslautern 2016), LMS Journal of Computation and Mathematics {\bf 19}\br(2016)\brr12--28
\bibitem[EJ18]{EJ18}
Elsenhans, A.-S.\ and Jahnel, J.: {\em Real and complex multiplication on
$K3$~surfaces
via period integration}, To appear in: Experimental Mathematics, {\tt doi: 10.1080/\discretionary{}{}{}10586458.\discretionary{}{}{}2022.2061649}
\bibitem[EJ20a]{EJ20a}
Elsenhans, A.-S.\ and Jahnel, J.: {\em Explicit families of\/
$K3$~surfaces having real multiplication,} To appear in: Michigan Mathematical Journal, {\tt doi: 10.1307/mmj/20205878}
\bibitem[EJ20b]{EJ20b}
Elsenhans, A.-S.\ and Jahnel, J.: {\em Computations with algebraic surfaces,} in: Mathematical software--ICMS 2020, Lecture Notes in Comput.\ Sci.~12097, {\em Springer,} Cham~2020, 87--93
\bibitem[EJ21]{EJ21}
Elsenhans, A.-S.\ and Jahnel, J.: {\em Frobenius trace distributions for
$K3$~surfaces,} {\tt arXiv:\discretionary{}{}{}2102.10620}
\bibitem[vGG]{vGG}
von zur Gathen, J.\ and Gerhard, J.: Modern computer algebra, {\em Cambridge University Press,} New York~1999
\bibitem[vG]{vG}
van Geemen, B.: {\em Some remarks on Brauer groups of
$K3$~surfaces,}
Adv.\ Math.\ {\bf 197}\br(2005)\brr222--247
\bibitem[GH]{GH}
Griffiths, P.\ and Harris, J.: Principles of algebraic geometry, {\em Wiley-Inter\-science,} New York~1978
\bibitem[Ha]{Ha}
Harvey, D.: {\em Computing zeta functions of arithmetic schemes,} Proc.\ Lond.\ Math.\ Soc.\ {\bf 111}\br(2015)\brr1379--1401
\bibitem[HS]{HS}
Harvey, D.\ and Sutherland, A.: {\em Computing Hasse-Witt matrices of hyperelliptic curves in average polynomial time, II.\ Frobenius distributions: Lang-Trotter and Sato-Tate conjectures,} in: Contemp.\ Math.~663, {\em AMS,} Providence~2016, 127--147
\bibitem[Ke]{Ke}
Kedlaya, K.: {\em Counting points on hyperelliptic curves using Monsky-Washnitzer cohomology,} Journal of the Ramanujan Mathematical Society {\bf 16}\br(2001)\brr323--338
\bibitem[Kn]{Kn}
Kneser, M.: Quadratische Formen, {\em Springer,} Berlin~2002
\bibitem[Mi]{Mi}
Milne, J.\,S.: {\em On a conjecture of Artin and Tate,} Ann.\ of Math.\ {\bf 102}\br(1975)\brr517--533
\bibitem[Ne]{Ne}
Neukirch, J.: Algebraic number theory, Grundlehren der Mathematischen Wis\-sen\-schaf\-ten 322, {\em Springer,} Berlin~1999
\bibitem[Og]{Og}
Ogus, A.: {\em Supersingular
$K3$
crystals,} in: Journ\'ees de G\'eom\'etrie Alg\'ebrique de Rennes (Rennes~1978) II, Ast\'erisque~64, {\em SMF,} Paris~1979, 3\--86
\bibitem[Pi]{Pi}
Pila, J.: {\em Frobenius maps of abelian varieties and finding roots of unity in finite fields,} Math.\ Comp.\ {\bf 55}\br(1990)\brr745--763
\bibitem[Sch]{Sch}
Schoof, R.: {\em Counting points on elliptic curves over finite fields,} J.\ Th\'eor.\ Nombres Bordeaux {\bf 7}\br(1995)\brr219--254
\bibitem[SGA4]{SGA4}
Artin, M., Grothendieck, A.\ et Verdier, J.-L.\ (avec la collaboration de Deligne, P.\ et Saint-Donat,~B.): Th\'eorie des topos et cohomologie \'etale des sch\'emas, S\'eminaire de G\'eom\'etrie Alg\'ebrique du Bois Marie
1963--1964 (SGA\,4), Lecture Notes in Math.~269, 270, 305, {\em Springer,} Berlin, Heidelberg, New York~1972--1973
\bibitem[SGA4\textonehalf]{SGA412}
Deligne, P.\ (avec la collaboration de Boutot, J.\ F., Grothendieck, A., Illusie, L.\ et Verdier, J.-L.): Cohomologie \smash{\'Etale,} S\'eminaire de G\'eom\'etrie Alg\'ebrique du Bois Marie
(SGA4\textonehalf),
Lecture Notes in Math.~569, {\em Springer,} Berlin, Heidelberg, New York~1977
\bibitem[SGA5]{SGA5}
Grothendieck, A.\ (avec la collaboration de Bucur, I., Houzel, C., Illusie, L.\ et Serre, \mbox{J.-P.}): Cohomologie
\mbox{$l$-adique}
et
Fonctions~$L$,
S\'eminaire de G\'eom\'etrie Alg\'ebrique du Bois Marie 1965--1966 (SGA\,5), Lecture Notes in Math.~589, {\em Springer,} Berlin, Heidelberg, New York~1977
\bibitem[Sk]{Sk}
Skorobogatov, A.\,N.: {\em Cohomology and the Brauer group of double covers,}
in: Brauer groups and obstruction problems: moduli spaces and arithmetic
(A.~Auel, B.~Hassett, A.~V\'arilly-Alvarado, and B.~Viray eds.), {\em Springer,} Cham~2017
\bibitem[Wi]{Wi}
Wilson, R.\,A.: The finite simple groups, Graduate Texts in Mathematics~251, {\em Springer,} London~2009
\bibitem[Yo]{Yo}
Yoshikawa, K.-I.: {\em Discriminant of certain
$K3$~surfaces,}
in: Representation theory and automorphic forms, Progr.\ Math.~255, {\em Birkh\"auser,} Boston~2008, 175--210
\end{thebibliography}
\end{document}